\documentclass[10pt]{amsart}
\usepackage{amsfonts,amssymb}
\RequirePackage{ifpdf}
\usepackage{amsmath}
\usepackage{color}
\usepackage{pstcol,pst-node}
\usepackage{graphics}
\usepackage{epic}
\usepackage{curves}

\def\appendix#1{
\addtocounter{section}{1} \setcounter{equation}{0}
\renewcommand{\thesection}{\Alph{section}}
\section*{Appendix \thesection\protect\indent\quad
#1}
}
\renewcommand{\theequation}{\thesection.\arabic{equation}}

\catcode`\@=11
\def\marginnote#1{}

\newcount\hour
\newcount\minute
\newtoks\amorpm
\hour=\time\divide\hour by60 \minute=\time{\multiply\hour by60
\global\advance\minute by-\hour}
\edef\standardtime{{\ifnum\hour<12 \global\amorpm={am}%
        \else\global\amorpm={pm}\advance\hour by-12 \fi
        \ifnum\hour=0 \hour=12 \fi
        \number\hour:\ifnum\minute<10 0\fi\number\minute\the\amorpm}}
\edef\militarytime{\number\hour:\ifnum\minute<100\fi\number\minute}

\def\draftlabel#1{{\@bsphack\if@filesw {\let\thepage\relax
      \xdef\@gtempa{\write\@auxout{\string
          \newlabel{#1}{{\@currentlabel}{\thepage}}}}}\@gtempa \if@nobreak
    \ifvmode\nobreak\fi\fi\fi\@esphack} \gdef\@eqnlabel{#1}}
    \def\@eqnlabel{}
\def\@vacuum{}
\def\draftmarginnote#1{\marginpar{\raggedright\scriptsize\tt#1}}

\def\draft{
  \oddsidemargin -.5truein
  \def\@oddfoot{\footnotesize \sl preliminary draft \hfil
    \rm\thepage\hfil\sl\today\quad\militarytime}
  \let\@evenfoot\@oddfoot \overfullrule 3pt
    \let\label=\draftlabel
    \let\marginnote=\draftmarginnote
  \def\@eqnnum{(\theequation)\rlap{\kern\marginparsep\tt\@eqnlabel}%
    \global\let\@eqnlabel\@vacuum}

  }

\def\eps{\varepsilon}
\def\be{\begin{equation}}
\def\ee{\end{equation}}
\def\bea{\begin{eqnarray}}
\def\eea{\end{eqnarray}}
\def\<{\langle}
\def\>{\rangle}

\def\nn{\nonumber}

\def\one#1{#1^{\raise5pt\hbox{$\scriptstyle\!\!\!\!1$}}\,{}}
\def\two#1{#1^{\raise5pt\hbox{$\scriptstyle\!\!\!\!2$}}\,{}}
\def\onetwo#1{#1^{\raise5pt\hbox{$\scriptstyle\!\!\!\!\!{12}$}}\,{}}

\newtheorem{theorem}{Theorem}[section]
\newtheorem{lm}[theorem]{Lemma}
\newtheorem{prop}[theorem]{Proposition}
\theoremstyle{definition}
\newtheorem{df}[theorem]{Definition}

\newtheorem{remark}[theorem]{Remark}

\newtheorem{cor}[theorem]{Corollary}
\theoremstyle{remark}

\allowdisplaybreaks

\begin{document}
\title[{\it Painlev\'e equations, Cherednik algebras and q-Askey scheme}]
{Confluences of the Painlev\'e equations, Cherednik algebras and q-Askey scheme.}
\author{Marta Mazzocco$^\dagger$}\thanks{$^\dagger$Department of Mathematical Sciences, Loughborough University, UK. Email: m.mazzocco@lboro.ac.uk, Phone: +44 1509 223187, 
Fax: +44 (0)1509 223969.}

\maketitle

\begin{abstract}
In this paper we produce seven new algebras as confluences of the Cherednik algebra of type  $\check{C_1}C_1$ and we characterise their spherical--sub-algebras. 

The limit of the spherical sub-algebra of the Cherednik algebra of type $\check{C_1}C_1$ is the monodromy manifold of the Painlev\'e VI equation \cite{O}. Here we prove that by considering the limits of the spherical sub-algebras of our new confluent algebras, one obtains the monodromy manifolds of all other Painlev\'e differential equations. Moreover, we introduce confluent versions of the Zhedanov algebra and prove that each of them (quotiented by their Casimir) is isomorphic to the corresponding spherical sub-algebra of our new confluent Cherednik algebras. We show that in the basic representation our confluent Zhedanov algebras act as symmetries of certain elements of the q-Askey scheme, thus setting a stepping stone towards the solution of the open problem of finding the corresponding quantum algebra for each element of the q-Askey scheme.

These results establish a new link between the theory of the Painlev\'e equations and the theory of the q-Askey scheme and shed light on the reasons behind the occurrence of special polynomials in the Painlev\'e theory.

\end{abstract}

\tableofcontents

\section{Introduction}
The relationship between the theory of the Painlev\'e equations and special or orthogonal polynomials is a very famous one and could be resumed by saying that thanks to the $\tau$--function structure of the Painlev\'e equations, some of their special solutions are related to special or orthogonal polynomials either directly, i.e. some rational solutions of the Painlev\'e equations are ratios of  special polynomials \cite{Yab,Vor,Ok1,Ok2,Ok3,NOOU,NY,Yam,Um,Cl}, or indirectly, i.e. some random matrix integrals which can be expressed by classical orthogonal polynomials have Fedholm determinants which can be expressed in terms of special solutions of the Painlev\'e equations \cite{TW,BI,FW,BD}.

It this paper we present a new relation between the theory of the Painlev\'e equations and q-polynomials belonging to the q-Askey scheme \cite{KLS}. This link does not rely on the $\tau$-function structure nor on choosing special solutions, it is indeed a much deeper and more conceptual relation that has allowed the author to discover seven new confluent Cherednik algebras and to prove several interesting results about them. 

Let us start from the Painlev\'e sixth equation \cite{fuchs, Sch, Gar1} which describes the monodromy preserving deformations of a rank $2$ Fuchsian system  with four simple poles $a_1,a_2,a_3$ and $\infty$. The solution of this Fuchsian system is in general a multi-valued analytic vector--function 
in the punctured Riemann sphere $\mathbb P^1\setminus\{a_1,a_2,a_3,\infty\}$ and its multivaluedness is 
described by the so-called monodromy group, i.e. a subgroup of $SL_2(\mathbb C)$ generated by the images $M_1,M_2,M_3$ of the generators of the fundamental group under the anti-homomorphism:
$$
\rho:\pi_1\left( \mathbb P^1\backslash\{a_1,a_2,a_3,\infty\},\lambda_0\right)\to SL_2(\mathbb C).
$$
The moduli  space ${\mathcal M}\slash{\Gamma}$ of monodromy representations $\rho$ up to Jordan equivalence, with prescribed local monodromy (i.e. prescribed conjugacy class for each $M_1,M_2,M_3$),  is realised as an affine cubic surface \cite{jimbo} which coincides with the center of the Cherednik algebra of type $\check{C_1}C_1$ for $q=1$ \cite{O,EOR}.  
In \cite{CM} this affine cubic surface was explicitly quantised leading to the Zhedanov algebra, which is to be isomorphic \cite{K1} to the spherical subalgebra of the Cherednik algebra of type $\check{C_1}C_1$ i.e. the algebra $\mathcal H$ generated by four elements $T_0,T_1,X,W$ with the following relations\footnote{Note that here $W=X^{-1}$. However we prefer here use the notation $W$ as in the confluent cases the generator $X$ may become singular.}
 \cite{Cher,Sa,NS,St}:
 \bea
\label{sahi1}
X W=W X=1,\\
\label{sahi2}
(T_1+ a b)(T_1+1)=0,\\
\label{sahi3}
(T_0+q^{-1} c d)(T_0+1)=0,\\
\label{sahi4}
(T_1 X +a)(T_1 X +b)=0,\\
\label{sahi5}
(qT_0 W+ c)(qT_0 W+d)=0,
\eea
where  $a,b,c,d,q\in\mathbb C^\star$.

In this paper, we produce seven new algebras as confluences of the Cherednik algebra of type  $\check{C_1}C_1$ in such a way that their spherical--sub-algebras tend in the limit to the monodromy manifolds of all other Painlev\'e differential equations:

\begin{df}\label{de:main}
Let $a,b,c,q\in\mathbb C^\star$. 
The confluent Cherednik algebras $\mathcal H_V,\mathcal H_{IV}$, $\mathcal H_{III}$,  $\mathcal H_{III^{D_7}}$,  $\mathcal H_{III^{D_8}}$, $\mathcal H_{II}$ are the 
unital associative $\mathbb C$--algebras generated by four elements $X,W,T_0,T_1$ satisfying the following relations respectively:
 \begin{itemize}
 \item $\mathcal H_V$:
 \bea
\label{sahi1-V}
X W= W X=1,\\
\label{sahi2-V}
(T_1+ a b)(T_1+1)=0,\\
\label{sahi3-V}
T_0(T_0+1)=0,\\
\label{sahi4-V}
(T_1 X +a)(T_1 X +b)=0,\\
\label{sahi6-V}
q T_0 W+ c = X (T_0+1),
\eea
\item $\mathcal H_{IV}$:
  \bea
\label{sahi1-IV}
X W= W X=0,\\
\label{sahi2-IV}
(T_1- b^2)(T_1+1)=0,\\
\label{sahi3-IV}
T_0(T_0+1)=0,\\
 \label{sahi7-IV}
T_1 X + b=W (T_1-b^2+ 1),\\
\label{sahi6-IV}
 q T_0 W + c=X (T_0+1),
\eea
\item $\mathcal H_{III}$:
  \bea
\label{sahi1-III}
X W= W X=1,\\
\label{sahi2-III}
(T_1+ a b)(T_1+1)=0,\\
\label{sahi3-III}
T_0^2=0,\\
\label{sahi4-III}
(T_1 X +a)(T_1 X +b)=0,\\
\label{sahi6-III}
q T_0 W-\sqrt{q}= X T_0.
\eea
 \item $\mathcal H_{III^{D_7}}$:
  \bea
\label{sahiPIIID71}
X W=W X=1,\\
\label{sahiPIIID72}
T_1( T_1+1)=0,\\
\label{sahiPIIID73}
T_0^2=0,\\
\label{sahiPIIID74}
T_1 X +a-W( T_1+1)   =0,\\
\label{sahiPIIID75}
q T_0 W-\sqrt{q}= X T_0,
\eea
 \item $\mathcal H_{III^{D_8}}$:
  \bea
\label{sahiPIIID81}
X W=W X=1,\\
\label{sahiPIIID82}
T_1( T_1+1)=0,\\
\label{sahiPIIID83}
T_0^2=0,\\
\label{sahiPIIID84}
T_1 X -W (T_1+1)  =0,\\
\label{sahiPIIID85}
q T_0 W-\sqrt{q}= X T_0.
\eea
\item $\mathcal H_{II}$:
  \bea
\label{sahi1-II}
X W= W X=0,\\
\label{sahi2-II}
(T_1-b^2)(T_1+1)=0,\\
\label{sahi3-II}
T_0(T_0+1)=0,\\
 \label{sahi7-II}
T_1 X =W (T_1-b^2+ 1),\\
\label{sahi6-II}
 q T_0 W -\sqrt{q}=X (T_0+1).
\eea
\end{itemize}
\end{df}
 
 It is worth noting that in the above definition, we have three different algebras corresponding to three different types of the third Painlev\'e equation \cite{sakai}.
 
The next set of results regards equivalent presentations for these confluent algebras.  For the 
Cherednik algebra of type $\check{C_1}C_1$ the following result is well known \cite{NS}:

\begin{lm}\cite{NS,O} The Cherednik algebra of type $\check{C_1}C_1$  is isomorphic as algebra to the unital associative $\mathbb C$--algebra with generators 
$V_0$, $\check{V_0}$, $V_1$, $\check{V_1}$ and relations:
\bea
\label{daha1}
(V_0-k_0)(V_0+k_0^{-1})=0\\
\label{daha2}
(V_1-k_1)(V_1+k_1^{-1})=0\\
\label{daha3}
(\check{V_0}-u_0)(\check{V_0}+u_0^{-1})=0\\
\label{daha4}
(\check{V_1}-u_1)(\check{V_1}+u_1^{-1})=0\\
\label{daha5}
\check{V_1}V_1V_0 \check{V_0}=q^{-1/2},
\eea
where  $k_0,k_1,u_0,u_1,q\in\mathbb C^\star$. The isomorphism is explicitly given by 
\be\label{eq:SO}
\phi(T_0)=k_0 {V_0},\quad \phi(T_1)= u_1\check{V_1},\quad \phi(X) = q^{1/2} V_0\check{V_0}, \quad \phi(W)= \check{V_1} V_1,
\ee
and for the parameters
\be\label{eq:SO1}
a=-\frac{u_1}{k_1},\quad b=k_1 u_1,\quad c=-q^{\frac{1}{2}} \frac{k_0}{u_0},\quad d=q^{1/2} u_0 k_0.
\ee
\end{lm}

In Theorem \ref{lm:sahi-others} we prove that all the confluent algebras $\mathcal H_V,\dots,\mathcal H_{II}$ are also isomorphic as algebras to some unital associative $\mathbb C$--algebras with generators $V_0$, $\check{V_0}$, $V_1$, $\check{V_1}$ and confluent analogues of the relations (\ref{daha1},\ref{daha2},\ref{daha3},\ref{daha4},\ref{daha5}). This allows to produce one further confluent algebra that was not visible in the previous presentation:

\begin{df}
Let $q\in\mathbb C^\star$. The confluent Cherednik algebra
$\mathcal H_{I}$ is the algebra 
with generators $V_0,\check{V_0},V_1,\check{V_1}$ and relations:
\bea
\label{daha1-P1}
V_0(V_0+1)=0\\
\label{daha2-P1}
V_1^2=0\\
\label{daha3-P1}
\check{V_0}(\check{V_0}+1)=0\\
\label{daha4-P1}
\check{V_1}(\check{V_1}+1)=0\\
\label{daha5-P1}
q^{1/2}\check{V_1}V_1V_0 = \check{V_0}+1,\\
\label{daha6-P1}
 V_0 \check{V_0}=0,
\eea
\end{df}

Next, we deal with the spherical sub--algebras of each confluent Cherednik algebra. We start by selecting a symmetriser $e$ and three special elements denoted by $X_1,X_2,X_3$ such that $[e,X_i]=0$, $i=1,2,3$, and such that $\hat X_i:= e X_i$, $i=1,2,3$, generate the spherical sub--algebras of each confluent Cherednik algebra. We prove that such elements $X_1,X_2,X_3$ satisfy a cubic relation (see Propositions \ref{pr:cubic-V}, \ref{pr:cubic-IV}, \ref{pr:cubic-III}, \ref{pr:cubic-IIID7}, \ref{pr:cubic-IIID8}, \ref{pr:cubic-II}, \ref{pr:cubic-I}) and that
 in the limit such cubic relations coincide with the monodromy manifolds of the corresponding Painlev\'e equations as defined in \cite{SvdP,vdP} (see Corollaries  \ref{co:cubic-V}, \ref{co:cubic-IV}, \ref{co:cubic-III}, \ref{co:cubic-IIID7}, \ref{co:cubic-IIID8}, \ref{co:cubic-II}, \ref{co:cubic-I}). In other words, one could say that the confluent Cherednik algebras introduced in this paper are such that their spherical sub--algebras produce a quantisation of the monodromy manifolds of the respective Painlev\'e equations. 
 
In order to link  the spherical sub--algebras $e {\mathcal H}_{V} e,\dots,e{\mathcal H}_{I} e$ to the q-Askey scheme, we first need to introduce the confluent versions of the Zhedanov algebra \cite{Zhe}:

\begin{df}\label{df-zhe}
The confluent Zhedanov algebras ${\mathcal Z}_{V}$, ${\mathcal Z}_{IV}$, ${\mathcal Z}_{III}$,   ${\mathcal Z}_{III}^{D_7}$, ${\mathcal Z}_{III}^{D_8}$, ${\mathcal Z}_{II}$, ${\mathcal Z}_{I}$ are the algebras generated by three elements $K_0$, $K_1$ and $K_2$  which satisfy the following relations:
\bea\label{zhe1-pv}
q^{\frac{1}{2}} K_0 K_1- q^{-\frac{1}{2}} K_1 K_0=K_2,\\
\label{zhe2-pv}
q^{\frac{1}{2}} K_1 K_2-q^{-\frac{1}{2}} K_2 K_1=B K_1+C_0 K_0+D_0,\\
\label{zhe3-pv}
q^{\frac{1}{2}}K_2 K_0-q^{-\frac{1}{2}}K_0 K_2=B K_0+D_1,
\eea
where the parameters $B,C_0,D_0,D_1$ are arbitrary non zero constants,  or have a specific form according to the following table:
\bea\label{eq:zhe-param}
&&
B=\left\{\begin{array}{cc}
\hbox{arbitrary} &\hbox{for } {\mathcal Z}_{V},{\mathcal Z}_{IV},\mathcal Z_{III}, \mathcal Z_{II}, \\
-\frac{(q-1)^2}{q}&\hbox{for } {\mathcal Z}_{III}^{D_7},{\mathcal Z}_{III}^{D_8},\\
-\frac{(q-1)^2}{\sqrt{q}}&\hbox{for } {\mathcal Z}_{I},\\
\end{array}\right.
\nn\\
&&
C_0=\left\{\begin{array}{cc}
\left(q-\frac{1}{q}\right)^2, &\hbox{for } {\mathcal Z}_{V},{\mathcal Z}_{III},{\mathcal Z}_{III}^{D_7},{\mathcal Z}_{III}^{D_8},\\
0,&\hbox{for } \mathcal Z_{IV}, \mathcal Z_{II},\mathcal Z_{I} \\
\end{array}\right.
\nn\\
&&
D_0=\left\{\begin{array}{cc}
\hbox{arbitrary}  &\hbox{for } {\mathcal Z}_{V},{\mathcal Z}_{IV},{\mathcal Z}_{III},{\mathcal Z}_{III}^{D_7}, \mathcal Z_{II}\\
0 &\hbox{for } \mathcal Z_{III}^{D_8},\mathcal Z_{I} \\
\end{array}\right.
\\
&&
D_1=\left\{\begin{array}{cc}
\hbox{arbitrary}, &\hbox{for } {\mathcal Z}_{V},{\mathcal Z}_{IV},\\
0,&\hbox{for } {\mathcal Z}_{III},{\mathcal Z}_{III}^{D_7},\mathcal Z_{III}^{D_8},\mathcal Z_{II}, \mathcal Z_{I}. \\
\end{array}\right.
\nn
\eea
\end{df}

\vskip 2mm
In Theorem \ref{th:zhe-sph}, we prove that the spherical sub-algebras of our confluent Cherednik algebras are isomorphic to the corresponding confluent Zhedanov algebras quotiented over their Casimir. 

In order to prove Theorem \ref{th:zhe-sph} we to provide an embedding for each (confluent) Cherednik algebras $\mathcal H, \mathcal H_V,\mathcal H_{IV},\mathcal H_{III},\mathcal H_{II}$ and $\mathcal H_{I}$ in $Mat(2,\mathbb T_q)$, where $T_q$ is the rank $2$ quantum torus, i.e. the ring of Laurent polynomials in $x,y,q$ with commutation relations $x y= q yx,\, q x= x q,\, q y= y q$ (see Theorems \ref{th-rep-qT},\ref{th-rep-qT-PV},\ref{th-rep-qT-PIV},\ref{th-rep-qT-PIII},\ref{th-rep-qT-PII}). These embeddings show that all of these algebras are Azumaya of rank $2$  (for $\mathcal H$, already appeared in \cite{O})\footnote{The author is grateful to Pavel Etingof for this observation.}.

Finally, we give a faithful representation of the confluent Zhedanov algebras  $\mathcal Z_V,\mathcal Z_{III}$, $\mathcal Z_{III}^{D_7},\mathcal Z_{III}^{D_8}$ on the space of symmetric Laurent polynomials 
and of the confluent Zhedanov algebras  $\mathcal H_{IV},\mathcal H_{II}$ and $\mathcal H_{I}$ on the space of polynomials and prove that specific elements of the q-Askey scheme arise as eigenvectors of one of the operators in such representations.\footnote{The basic representation for the algebras  ${\mathcal Z}_{V}$, 
${\mathcal Z}_{III}$,   ${\mathcal Z}_{III}^{D_7}$, ${\mathcal Z}_{III}^{D_8}$ can be found in \cite{MS}.} There results are schematically  resumed in  figure 1.

  \begin{figure}[h]
\begin{pspicture}(-2,6)(2,-4)
 \rput(0,5){$\begin{array}{c}\hbox{Askey--Wilson}\\
 \mathcal Z_{VI}\\ \end{array}$}
\psline{->}(0,4.5)(0,3.6)\rput(0,3){$\begin{array}{c}
 \mathcal Z_{V}\\\qquad\qquad \qquad\hbox{Big q-Jacobi }\qquad \hbox{Continuous dual $q$--Hahn}\\\end{array}$} 
  \psline{->}(-1.5,2.5)(-2.5,1.5)
  \rput(-2.5,1){$\begin{array}{c}\hbox{Big $q$--Laguerre}\\
 \mathcal Z_{IV}\\ \end{array}$}
 \psline{->}(1.5,2.5)(2.5,1.5)
  \rput(2.5,1){$\begin{array}{c}\hbox{Al Salam--Chihara}\\
 \mathcal Z_{III}\\ \end{array}$}
   \psline{->}(-2.5,0.5)(-2.5,-.5)
  \rput(-2.5,-1){$\begin{array}{c}\hbox{Little $q$--Laguerre}\\
 \mathcal Z_{II}\\ \end{array}$}
   \psline{->}(-2.5,-1.5)(-2.5,-2.5)
  \rput(-2.5,-3){$\begin{array}{c}\hbox{Little $q$--Laguerre with $a=0$}\\
 \mathcal Z_{I}\\ \end{array}$}
   \psline{->}(2.5,0.5)(2.5,-.5) 
   \rput(2.5,-1){$\begin{array}{c}\hbox{\qquad Continuous big $q$--Hermite}\\
 \mathcal Z_{III}^{D_7}\\ \end{array}$}
  \psline{->}(2.5,-1.6)(2.5,-2.5) 
   \rput(2.5,-3){$\begin{array}{c}\hbox{Continuous $q$--Hermite}\\
 \mathcal Z_{III}^{D_8}\\ \end{array}$}
   \end{pspicture}
   \caption{The confluence scheme for the Zhedanov algebras and the polynomials in the q-Askey scheme. Note that each confluent Zhedanov algebra is labelled by the Painlev\'e equation whose the monodromy manifold arises as its semiclassical limit.} \end{figure}
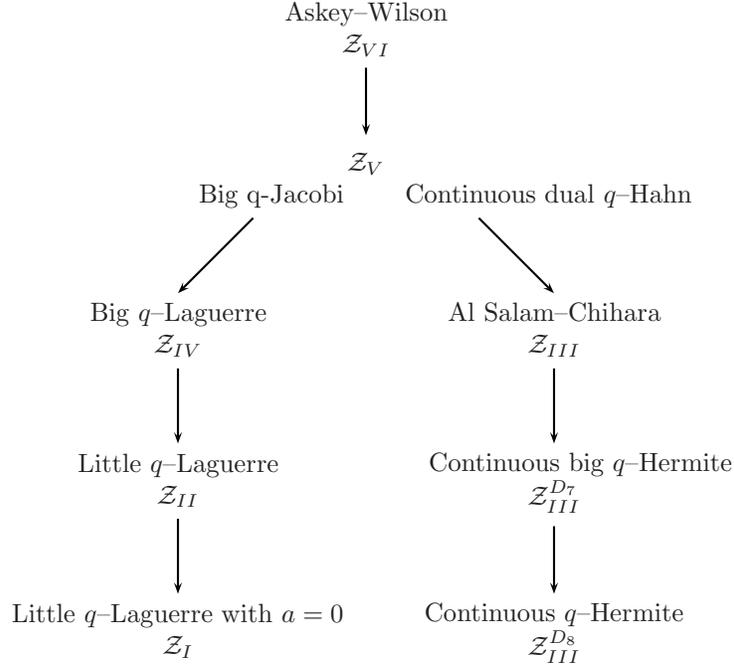

Note that for ${\mathcal Z}_{V}$ we have two different faithful representations corresponding to the continuous dual $q$--Hahn polynomials and to the big 
$q$--Jacobi polynomials, which is due to an algebra automorphism 
of $\mathcal H$ which under confluence produces an algebra isomorphic to $\mathcal H_V$ as described in sub--section \ref{der} and in Lemma \ref{auto-Zhe}. This algebra isomorphism reflects the duality between continuous dual $q$--Hahn polynomials and big 
$q$--Jacobi polynomials \cite{K4}. The author conjectures that all families of polynomials in the left side of the q-Askey scheme can be obtained with a similar construction from confluences of the algebra automorphisms 
of $\mathcal H$.

\begin{remark}
In \cite{Cher1}, Cherednik introduced a Whittaker limit on difference spherical functions and correspondingly introduced a notion of {\it nil-DAHA} for the root system $A_n$. The confluent Cherednik algebras defined in the current paper are a concatenation of Whittaker-type limits in the case of root system $C_1$.\footnote{The author is grateful to P. Etingof for pointing this out to her.} \end{remark}

This paper  is organised as follows:  In Section \ref{se:background}, we recall some background material on the theory of the Cherednik algebra of type  $\check{C_1}C_1$.  In Section \ref{se:first-results}, we explain how to derive our confluent Cherednik algebras and 
give some equivalent presentations for them.  In Section \ref{se:spherical-others}, we discuss the spherical sub--algebras and show that the generators satisfy a cubic relation which in the limit coincides with the monodromy manifolds of the corresponding Painlev\'e equations. In Section \ref{se:embedding} we provide  an embedding for each (confluent) Cherednik algebras $\mathcal H, \mathcal H_V,\mathcal H_{IV},\mathcal H_{III},\mathcal H_{II}$ and $\mathcal H_{I}$ in $Mat(2,\mathbb T_q)$, where $\mathbb T_q$ is the rank $2$ quantum torus. In Section \ref{se:zhedanov-others}, we prove that each spherical sub-algebra is isomorphic to the corresponding confluent Zhedanov algebra quotiented by its Casimir and show that the latter act as symmetries of some elements of the q-Askey scheme. In the appendix we list the Painlev\'e monodromy manifolds.

\section{Notation and background on the Cherednik algebra of type $\check{C_1}C_1$}\label{se:background}

In this section we recall some background material on the theory of the Cherednik algebra of type $\check{C_1}C_1$ and  the relation due to Koornwinder \cite{K1,K2} between its spherical sub--algebra and Askey--Wilson polynomials. Throughout this section we assume $q^m\neq 1$, for all $m\in\mathbb Z$. 

\subsection{Automorphisms of the Cherednik algebra  of type $\check{C_1}C_1$}\label{suse:aut}

The automorphisms of  the Cherednik algebra  of type $\check{C_1}C_1$ were studied in \cite{NS,St}. Here we list the ones that will be used in this paper:

\begin{prop}
The following transformations are automorphisms of the Cherednik algebra of type $\check{C_1}C_1$:
\bea
&&
\beta: \left\{\begin{array}{lcl}
T_0&\to&-\frac{q}{c} W T_0-\left(1+\frac{d}{c}\right),\\
T_1&\to&T_1,\\
X&\to&X,\\
a &\to& a,\\
b &\to& b,\\
c &\to& \frac{q}{c},\\
d &\to& d,\\
\end{array}
\right.\nn\\
&&
\gamma: \left\{\begin{array}{lcl}
T_0&\to&b T_1^{-1} W,\\
T_1&\to&T_1,\\
X&\to&\sqrt{\frac{a b c d}{q}}T_1^{-1} W T_0^{-1} X,\\
a &\to& \sqrt{\frac{ a b c d}{q}},\\
b &\to& -\sqrt{\frac{q a b}{ c d}},\\
c &\to&-\sqrt{\frac{q b c}{a d}},\\
d &\to&\sqrt{\frac{q b d}{a c}}.\\
\end{array}
\right.\eea
\end{prop}

\subsection{Zhedanov algebra and Askey Wilson polynomials}\label{se:basic}

Let us remind the definition of Zhedanov algebra:

\begin{df}\cite{Zhe}
Let $B,C_0,C_1,D_0,D_1\in\mathbb C$. The Zhedanov algebra ${\mathcal Z}$ is the algebra generated by three elements $K_0$, $K_1$ and $K_2$  which satisfy the following relations:
\bea\label{zhe1}
q^{\frac{1}{2}} K_0 K_1- q^{-\frac{1}{2}} K_1 K_0=K_2,\\
\label{zhe2}
q^{\frac{1}{2}} K_1 K_2-q^{-\frac{1}{2}} K_2 K_1=B K_1+C_0 K_0+D_0,\\
\label{zhe3}
q^{\frac{1}{2}}K_2 K_0-q^{-\frac{1}{2}}K_0 K_2=B K_0+C_1 K_1+D_1.
\eea
\end{df}
Note that this algebra admits the following Casimir
\bea\label{eq:casimir-zhe}
&&
\mathcal C:= q^{-\frac{1}{2}}(1-q^2)K_0 K_1 K_2+q K_2^2+B(K_0 K_1+K_1 K_0)+qC_0K_0^2+\nn\\
&&
+\frac{C_1}{q} K_1^2+(1+q)D_0K_0+(1+\frac{1}{q})D_1 K_1.
\eea
The Zhedanov algebra depends on $5$ parameters, but we can choose two of them, for example $C_1$ and $C_0$ by rescaling the generators. The quotient algebra
$\mathcal Z(\mathcal C_0)={\mathcal Z}\slash\langle \mathcal C= \mathcal C_0\rangle$ therefore depends on $4$ independent parameters (one of them being the value of 
$\mathcal C_0$).

In \cite{K1} Koornwinder defined an embedding of the central extension of the Zhedanov algebra $\mathcal Z(\mathcal C_0)$, into the Cherednik algebra $\mathcal H$ of type $\check{C_1}C_1$. This result was then generalised to the universal Askey--Wilson algebra defined in \cite{Ter} in 
\cite{Ter1}. 

Without going into too much detail, let us recall the main ingredients of Koornwinder embedding. Let us express the Zhedanov algebra structure constants by the parameters $a,b,c,d$:
\bea\label{eq:zhe-par-abcd}
&&
B = \frac{(q-1)^2}{q}\left(\left(1+\frac{a b}{q}\right) \left(\frac{d}{c}+1\right)c + 
\left(\frac{b}{a}+1\right) \left(1+\frac{c d}{q}\right)a\right),\nn\\
&&
C_0=\left(q-\frac{1}{q}\right)^2,\qquad C_1=\frac{a b c d}{q}\left(q-\frac{1}{q}\right)^2 \nn \\
&&
D_0=-\frac{(q+1)(q-1)^2}{q} \left(\left( \frac{b}{a}+1\right) \left(\frac{d}{c}+1\right)\frac{a c}{q} + 
\left(1+\frac{a b}{q}\right) \left(1+\frac{c d}{q}\right)\right),\nn\\
&&
D_1=-\frac{(q+1)(q-1)^2}{q^2} \left(\left( \frac{b}{a}+1\right) \left(1+\frac{a b}{q}\right)a  c d+ 
\left(\frac{d}{c}+1\right) \left(1+\frac{c d}{q}\right)a b c \right),\nn
\eea
then $\mathcal Z(\mathcal C_0)$ admits the following representation on the space $\mathcal L_{sym}$ of symmetric Laurent polynomials \cite{K1,K3}:
\bea
\label{eq:basic4}
&\qquad (K_1 f)[x]:=&( x+\frac{1}{x}) f[x],\\
\label{eq:basic5}
&\qquad(K_0 f)[x]:=& \frac{(1-a x)(1-b x)(1-c x)(1-d x)}{(1-x^2)(1-q x^2)}(f[qx]-f[x])+\nn\\
&&+\frac{(a-x)(b-x)(c-x)(d-x)}{(1-x^2)(q-x^2)}(f[q^{-1}x]-f[x])+\\
&&+(1+\frac{a b c d}{q}) f[x].\nn
\eea
We don't list here the representation for the operator $K_2$ as this can be obtained  by applying relation (\ref{zhe1}). 

The Askey Wilson  polynomials (we write them here in monic form like in \cite{K1}):
$$
P_n(x;a,b,c,d) := \frac{(ab,ac,ad;q)_n}{a^n(ab c d q^{n-1};q)_n}    {}_4\phi_3\left(\begin{array}{cc}
q^{-n},q^{n-1}abcd,ax,a x^{-1}\\
ab,ac,ad
\end{array};q,q
\right),
$$
are eigenfunctions of the $K_0$ operator:
$$
K_0 P_n =(q^{-n}+abcd q^{n-1}) P_n.
$$

The reduction from  the space $\mathcal L$ of  Laurent polynomials to  the space $\mathcal L_{sym}$ of symmetric Laurent polynomials is due to the action of the symmetriser of $\mathcal H$:
\be\label{df:e}
e:= \frac{1+i\sqrt{ab}\check{V_1}}{1-a b}
\ee
which allowed Koornwinder to establish the isomorphism between a central extension $\widetilde{\mathcal Z}(\mathcal C_0)$ and the spherical sub-algebra $e\mathcal He$  of $\mathcal H$. We discuss this result and the link with the PVI monodromy manifold in the next subsection.

\subsection{An important cubic relation and the spherical sub-algebra $e\mathcal H e$}\label{suse:cubic}

The following result plays an important role in this paper:

\begin{prop} \label{lm:spherical}
The following three elements:
\bea\label{spherical}
\noindent& X_1 =X+W,\qquad 
X_2 =\sqrt{\frac{a b c d}{q}}T_0^{-1}T_1^{-1} +\sqrt{\frac{q}{a b c d}} T_1T_0,\\
& X_3 =q \sqrt{a b c d} \left(\frac{1}{a b c d}W T_0 T_1 +\frac{1}{a b c d}X T_1 T_0+
\left(\frac{1}{a}+\frac{1}{b}\right)\frac{1}{c d} T_0+\left(\frac{1}{c}+\frac{1}{d}\right)\frac{1}{a b q} T_1
\right),
\nn\eea
satisfy the quantum commutation relations:
\bea\label{eq:skein}
&&q^{\frac{1}{2}} X_2 X_1-q^{-\frac{1}{2}} X_1 X_2 =  \left(q-\frac{1}{q}\right)X_3 -\nn\\
&&\qquad-  \left(q^{\frac{1}{2}} -q^{-\frac{1}{2}} \right)\left(\overline k_0\overline k_1- \overline u_0\left(\frac{i}{\sqrt{q a b}} T_1 + i \sqrt{q a b} T_1^{-1}\right)\right),\\
&&q^{\frac{1}{2}} X_3 X_2-q^{-\frac{1}{2}} X_2 X_3 =  \left(q-\frac{1}{q}\right)X_1 -\nn\\
&&\qquad-  \left(q^{\frac{1}{2}} -q^{-\frac{1}{2}} \right)\left(\overline k_0 \overline u_0-\overline k_1\left(\frac{i}{\sqrt{q a b}} T_1 + i \sqrt{q a b} T_1^{-1}\right)\right),\\
&&q^{\frac{1}{2}} X_1 X_3-q^{-\frac{1}{2}} X_3 X_1 =  \left(q-\frac{1}{q}\right)X_2 -\nn\\
&&\qquad- \left(q^{\frac{1}{2}} -q^{-\frac{1}{2}} \right)\left(\overline k_1 \overline u_0-\overline k_0\left(\frac{i}{\sqrt{q a b}} T_1 + i \sqrt{q a b} T_1^{-1}\right)\right),
\eea
and the quantum cubic relation:
\begin{eqnarray}\label{eq:skein-cubic}
q^{\frac{1}{2}}X_2 X_1 X_3-q X_2^2-\frac{1}{q} X_1^2-q X_3^2+ \sqrt{q}\left(\overline k_1 \overline u_0-\overline k_0 
\left(\frac{i}{\sqrt{q a b}} T_1 + i \sqrt{q a b} T_1^{-1}\right)\right) X_2+\nn\\
+ \frac{1}{\sqrt{q}}\left(\overline u_0\overline k_0-
\overline k_1\left(\frac{i}{\sqrt{q a b}} T_1 + i \sqrt{q a b} T_1^{-1}\right) \right)X_1+\nn\\
+\sqrt{q}\left(\overline k_0 \overline k_1-
\overline u_0\left(\frac{i}{\sqrt{q a b}} T_1 + i \sqrt{q a b} T_1^{-1}\right)\right) X_3+\\
+\overline k_0^2+\overline k_1^2+\overline u_0^2-\overline u_1^2+2\left(q+\frac{1}{q}\right)+ 
(\frac{q+1}{\sqrt{q}} \overline u_1-\overline k_0\overline k_1\overline u_0 )\left(\frac{i}{\sqrt{q a b}} T_1 + i \sqrt{q a b} T_1^{-1}\right),\nn
\end{eqnarray}
where:
\bea\label{eq:bar-parameters}
&&
\overline u_0= u_0-\frac{1}{u_0}= i\sqrt{\frac{d}{c}}+i\sqrt{\frac{c}{d}},\quad \overline k_0= k_0-\frac{1}{k_0}=i\sqrt{\frac{c d}{q}}+i\sqrt{\frac{q}{c d}},\\ 
&&
\overline u_1=u_1-\frac{1}{u_1}=i\sqrt{a b} +i\frac{1}{\sqrt{a b}},\quad \overline k_1= k_1-\frac{1}{k_1}=i\sqrt{\frac{b}{a}}+i\sqrt{\frac{a}{b}}.\nn
\eea
\end{prop}

\proof
Relations (\ref{eq:skein}) are proved in \cite{IT} {\rm (see also \cite{EE})}. The only proof the author could find of this cubic relation (\ref{eq:skein-cubic}) is in \cite{EE} for the case $\overline u_0=\overline u_1=\overline k_0=\overline k_1=0$. To prove it in the generic case, we use the following substitutions which can be easily be derived from (\ref{sahi2},\ref{sahi3},\ref{sahi4}) and (\ref{sahi5}):
\begin{eqnarray}\label{q-alg-ordering}
T_1^2=-(a b+1) T_1-a b, \qquad T_0^2=-(c d/q +1)T_0-\frac{cd}{q}\nn\\
T_1 X=W T_1+(1+ a b )W-(a+b),\\
 T_1 W= X T_1-(1+a b) W + a+b,\nn\\
T_0 X= q W T_0+c+d-(1+ cd/q)X,\nn\\
T_0 W =\frac{1}{q} X T_0+\frac{1}{q}(1+c d/q)X-\frac{1}{q}(c+d).\nn
\end{eqnarray}
Using these relations recursively, it is a straightforward computation to eliminate all higher powers in $T_0$ and $T_1$ and to bring all $X$ and $W$ to the left and eventually obtain the statement.
\endproof

The following lemma characterises the spherical--subalgebra $e \mathcal H e$:

\begin{cor}\cite{K2,EE,IT}
The elements $\hat X_i= e X_i e$, $i=1,2,3$, where $X_1,X_2,X_3$ are defined by (\ref{spherical}), generate the spherical sub-algebra $e\mathcal He$, they satisfy the quantum commutation relations:
\bea
\label{eq:skein-hat}
&&q^{\frac{1}{2}}\hat X_2 \hat X_1-q^{-\frac{1}{2}} \hat X_1 \hat X_2 =  \left(q-\frac{1}{q}\right)\hat X_3 - 
 \left(q^{\frac{1}{2}} -q^{-\frac{1}{2}} \right)\omega_3 e,\\
&&q^{\frac{1}{2}} \hat X_3 \hat X_2-q^{-\frac{1}{2}} \hat X_2 \hat X_3 =  \left(q-\frac{1}{q}\right)\hat X_1 - 
 \left(q^{\frac{1}{2}} -q^{-\frac{1}{2}} \right)\omega_1 e,\\
&&q^{\frac{1}{2}} \hat X_1 \hat X_3-q^{-\frac{1}{2}} \hat X_3 \hat X_1 =  \left(q-\frac{1}{q}\right)\hat X_2 -
 \left(q^{\frac{1}{2}} -q^{-\frac{1}{2}} \right)\omega_2 e,
\eea
and the following cubic relation:
\be\label{eq:cubic-hat}
q^{\frac{1}{2}} \hat X_2 \hat X_1\hat X_3 -
q \hat X_2^2-q^{-1}  \hat X_1^2- q \hat X_3^2+\sqrt{q} \omega_2 \hat X_2+ \frac{1}{\sqrt{q}} \omega_1 \hat X_1+\sqrt{q} \omega_3 \hat X_3-\omega_4 e=0.
\ee
where
\begin{eqnarray}\label{omegaPVI}
\omega_1=\overline u_0\overline k_0+ i
\overline k_1\left(\sqrt{\frac{a b}{q}}+\sqrt{\frac{q}{a b}}\right),\nn\\
\omega_2=\overline k_1\overline u_0+i \overline k_0\left(\sqrt{\frac{a b}{q}}+\sqrt{\frac{q}{a b}}\right),\\
\omega_3=\overline k_0 \overline k_1+
i \overline u_0 \left(\sqrt{\frac{a b}{q}}+\sqrt{\frac{q}{a b}}\right),\nn\\
\omega_4=\overline k_0^2+\overline k_1^2+\overline u_0^2-\left(\sqrt{\frac{a b}{q}}+\sqrt{\frac{q}{a b}}\right)^2 -
i\overline k_0 \overline k_1\overline u_0\left(\sqrt{\frac{a b}{q}}+\sqrt{\frac{q}{a b}}\right)+\frac{(1+q)^2}{q},\nn
\end{eqnarray}
where $\overline k_0, \overline k_1, \overline u_0, \overline u_1$ are defined in (\ref{eq:bar-parameters}).\end{cor}

\begin{lm}\cite{O}
In the limit $q\to 1$, 
$X_1, X_2, X_3$ satisfy the following cubic relation:
\bea\label{cubic}
&&
X_1  X_2 X_3 - X_1^2- X_2^2-X_3^2+(\overline u_0 \overline k_0+\overline u_1 \overline k_1)X_1+ (\overline k_1 \overline u_0+ \overline k_0 \overline u_1)  X_2+\nn\\
&&\quad
+(\overline k_0 \overline k_1+\overline u_0\overline u_1)  X_3 + \overline k_0^2+ \overline k_1^2+ \overline u_0^2+ \overline u_1^2-  \overline k_0 \overline k_1  \overline u_0  \overline u_1 +4=0.\nn
\eea
\end{lm}

\begin{remark}
This cubic is also known as the monodromy manifold of the sixth Painlev\'e equation (see Appendix A). In this paper we will obtain similar cubic relations for the spherical subalgebras of each confluent Cherednik algebra and we will show  that in the classical limit each of these cubic relations produces the monodromy manifold of the corresponding Painlev\'e equation.
\end{remark}

\section{Derivation and first properties of the confluent Cherednik algebras}\label{se:first-results}

The procedure to derive the confluent Cherednick algebras given in Definition \ref{de:main} can be resumed in two main roads, i.e.  $\mathcal H_{VI}\to \mathcal H_{V}\to \mathcal H_{III}\to \mathcal H_{III}^{D_7}\to \mathcal H_{III}^{D_8}$ and $\mathcal H_{VI}\to \mathcal H_{V}\to \mathcal H_{IV}\to \mathcal H_{II}\to \mathcal H_{I}$. In the first road we make $T_0$ more and more singular, then we make $T_1$ more and more singular, while keeping $X$ and $W$ invertible. In the second road we make $X$ and $W$ singular, then $T_0$ more and more singular. To reach $ \mathcal H_{I}$ however, we need to swap to the Noumi--Stokman presentation of the algebra $\mathcal H_{VI}$ otherwise we end up by loosing a generator. For this reason in this section we start by giving the Noumi--Stokman presentation for the algebras $\mathcal H_{VI}, \mathcal H_{V}, \mathcal H_{IV}, \mathcal H_{III} \mathcal H_{II}$, then we explain our derivation procedure more in detail.

\begin{remark}
As pointed out to the author by T. Koornwinder, there always is a degree of arbitrariness in such a confluence procedure. However there are two very strong mechanisms to remove such arbitrariness: the first one is that many confluences lead to algebras with too many relations and therefore are to be discarded. The second is that many different confluences give rise to algebras which are equivalent by an algebra isomorphism as in the case of $\mathcal H_V$. In this paper we impose a specific degeneration scheme for the cubic relations satisfied by the generators of the spherical sub--algebras such that in the limit they give rise to the Poisson relations on the monodromy manifolds of the Painlev\'e equations (see Section \ref{se:spherical-others}). 
\end{remark}

\begin{theorem}\label{lm:sahi-others}
Let $k_1,u_0,u_1,q\in\mathbb C^\star$. 
The confluent Cherednik algebras $\mathcal H_V,\mathcal H_{IV}$, $\mathcal H_{III},\mathcal H_{II}$ are  is isomorphic as algebras to the associative $\mathbb C$--algebra with generators $V_0,\check{V_0},V_1,\check{V_1}$ and relations:
 \begin{itemize}
 \item $\mathcal H_V$:
\bea
\label{dahaV1}
V_0^2+ V_0=0,\\
\label{dahaV2}
(V_1-k_1)(V_1+k_1^{-1})=0,\\
\label{dahaV3}
\check{V_0}^2+u_0^{-1}\check{ V_0}=0,\\
\label{dahaV4}
(\check{V_1}-u_1)(\check{V_1}+u_1^{-1})=0,\\
\label{dahaV5}
q^{1/2}\check{V_1}V_1V_0 =\check{V_0}+u_0^{-1},\\
\label{dahaV6}
q^{1/2}\check{V_0}\check{V_1}V_1=V_0+1.
\eea
\item $\mathcal H_{IV}$:
\bea
\label{dahaIV1}
V_0^2+ V_0=0,\\
\label{dahaIV2}
V_1^2+V_1=0,\\
\label{dahaIV3}
\check{V_0}^2+\frac{1}{u_0}\check{ V_0}=0,\\
\label{dahaIV4}
(\check{V_1}-u_1)(\check{V_1}+u_1^{-1})=0,\\
\label{dahaIV5}
q^{1/2}\check{V_1}V_1V_0 =\check{V_0}+u_0^{-1},\\
\label{dahaIV6}
\check{V_0}\check{V_1}V_1=0,\\
\label{dahaIV7}
V_0\check{V_0}=0.
\eea
\item $\mathcal H_{III}$:
\bea
\label{dahaIII1}
V_0^2=0,\\
\label{dahaIII2}
(V_1-k_1)(V_1+k_1^{-1})=0,\\
\label{dahaIII3}
\check{V_0}^2+\check{ V_0}=0,\\
\label{dahaIII4}
(\check{V_1}-u_1)(\check{V_1}+u_1^{-1})=0,\\
\label{dahaIII5}
q^{1/2}\check{V_1}V_1V_0 =\check{V_0}+1,\\
\label{dahaIII6}
q^{1/2}\check{V_0}\check{V_1}V_1=V_0.
\eea
\item $\mathcal H_{II}$:
\bea
\label{dahaII1}
V_0^2=0,\\
\label{dahaII2}
V_1^2+V_1=0,\\
\label{dahaII3}
\check{V_0}^2+\check{ V_0}=0,\\
\label{daha-lim4-3}
(\check{V_1}-u_1)(\check{V_1}+u_1^{-1})=0,\\
\label{dahaII4}
q^{1/2}\check{V_1}V_1V_0 =\check{V_0}+1,\\
\label{dahaII5}
\check{V_0}\check{V_1}V_1=0,\\
\label{dahaII6}
V_0\check{V_0}=0.
\eea
\end{itemize}
\end{theorem}

\proof
To avoid repeating the proof three times we use the following conventions:
$$\begin{array}{l}
u_0=\sqrt{q}\hbox{ for } \mathcal H_{III} \\
\\
u_0=1\hbox{ for } \mathcal H_{II} \\
\\
k_1=1\hbox{ for } \mathcal H_{IV}, \mathcal H_{II}.\\
\end{array}
$$
It is enough to give relations between the generators $V_0,V_1,\check{V_0},\check{V_1}$ and $X,W, T_0,T_1$:
\bea\label{eq:Sa-NS}&&
T_0= {V_0},\quad T_1= u_1\check{V_1},\quad  W= \check{V_1} V_1,\nn\\
&&
X = \left\{\begin{array}{lc}(V_1+k_1^{-1}-k_1)(\check{V_1}+u_1^{-1}-u_1),&\hbox{ for $\mathcal H_V$ and $\mathcal H_{III}$} \\
(V_1+1)(\check{V_1}+u_1^{-1}-u_1),&\hbox{ for $\mathcal H_{IV}$,} \\
\end{array}\right.
\eea
and for the parameters:
\be\label{eq:Sa-NS1}
a=-\frac{u_1}{k_1},
\quad b=k_1 u_1,\quad c=-q^{\frac{1}{2}} \frac{1}{u_0}.
\ee
Viceversa:
\be\label{eq:SO3}
{V_0}= T_0 ,\quad \check{V_1}=\frac{1}{u_1} T_1,\quad 
\check{V_0}= {q^{1/2}} W T_0-\frac{1}{u_0}, \quad  {V_1} =u_1 T_1^{-1}W,
\ee
where
\be\label{eq:inverse}
T_1^{-1}= -\frac{1}{a b} T_1 -(1+\frac{1}{a b}).
\ee
\endproof

\subsection{Derivation of $\mathcal H_{V}$}\label{der}

The derivation procedure is more transparent in Noumi--Stokman representation. Start from $\mathcal H$ and choose to rescale $V_0$ and $\check{V_0}$. Then (\ref{daha5}) will become singular, so before taking the limit we replace it by:
$$
\sqrt{q}\, \check{V_1}V_1V_0=\check{V_0}-\overline u_0,\qquad 
\sqrt{q}\, \check{V_0}\check{V_1}V_1={V_0}-\overline k_0.
$$
Now rescale: $V_0\to\frac{1}{\eps} V_0$, $\check{V_0}\to\frac{1}{\eps}\check{V_0}$, $k_0\to\eps$, and $u_0\to\eps u_0$. Then the defining relations (\ref{daha1},\ref{daha3},\ref{daha5}) become
\bea
\frac{1}{\eps^2}(V_0- \eps^2)(V_0+k_0^{-1})=0,\qquad
\frac{1}{\eps^2}(\check{V_0}-\eps^2 u_0)(\check{V_0}+u_0^{-1})=0,\nn\\
\frac{1}{\eps}\sqrt{q} \check{V_1}V_1V_0=\frac{1}{\eps}\check{V_0}+\frac{1}{\eps}u_0^{-1},\qquad 
\frac{1}{\eps}\sqrt{q} \check{V_0}\check{V_1}V_1=\frac{1}{\eps}{V_0}+\frac{1}{\eps}k_0^{-1}.\nn
\eea
and in the limit $\eps\to 0$ we obtain $\mathcal H_{V}$. Observe that the new $V_0$ and $\check{V_0}$ are no longer invertible and $q$ has not been rescaled - differently from the case of rational Cherednik algerbas \cite{BEG,G}.

We can derive $\mathcal H_{V}$ also in another way: choose to rescale $V_1$ and $\check{V_0}$. Then (\ref{daha5}) will become singular and needs to be replaced by:
$$
\sqrt{q}\, V_0 \check{V_0}\check{V_1}=\check{V_1}-\overline k_1\qquad 
\sqrt{q}\, \check{V_1}V_1V_0=\check{V_0}-\overline u_0,.
$$
Now rescale: $V_1\to\frac{1}{\eps} V_1$, $\check{V_0}\to\frac{1}{\eps}\check{V_0}$, $k_1\to-\frac{1}{\eps}$, and $u_0\to\eps u_0$. Then the defining relations (\ref{daha2},\ref{daha3},\ref{daha5}) become
\bea
\frac{1}{\eps^2}(V_1+ 1)(V_1-\eps^2)=0,\nn\\
\frac{1}{\eps^2}(\check{V_0}-\eps^2 u_0)(\check{V_0}+u_0^{-1})=0,\nn\\
\frac{1}{\eps}\sqrt{q}\, V_0 \check{V_0}\check{V_1}=\frac{1}{\eps}{V_1}+\frac{1}{\eps} \qquad 
\frac{1}{\eps}\sqrt{q}\, \check{V_1}V_1V_0=\frac{1}{\eps}\check{V_0}+\frac{1}{\eps}u_0^{-1}\nn.
\eea
By taking the limit $\eps\to 0$ we obtain the following algebra $\mathcal H_{V}^\gamma$:
\bea
\label{dahaV1gamma}
(V_0-k_0)(V_0+k_0^{-1})=0\\
\label{dahaV2gamma}
(V_1+1)V_1=0\\
\label{dahaV3gamma}
\check{V_0}^2+u_0^{-1}\check{ V_0}=0\\
\label{dahaV4gamma}
(\check{V_1}-u_1)(\check{V_1}+u_1^{-1})=0\\
\label{dahaV5gamma}
q^{1/2}\check{V_1}V_1V_0 =\check{V_0}+u_0^{-1}\\
\label{dahaV6gamma}
q^{1/2}V_0 \check{V_0}\check{V_1}={V_1}+1
\eea
This  algebra is the image of $\mathcal H_{V}$  in the (\ref{dahaV1}-\ref{dahaV6}) presentation by the following isomorphism $\gamma$:
\be\label{eq:gamma-V}
\gamma(\check{V_1},V_1,V_0,\check V_0)= (\check{V_1},V_1V_0 V_1^{-1},V_1,\check{V_0}),\quad
\gamma(u_1,k_1,k_0,u_0) =  (u_1,k_0^{-1},k_1,u_0),
\ee
where we pick $k_0=1$. In Sahi presentation the algebra $\mathcal H_{V}^\gamma$ takes the following form:
 \bea
\label{sahi1-Vgamma}
X W= W X=0,\\
\label{sahi2-Vgamma}
(T_1-a^2)(T_1+1)=0,\\
\label{sahi3-Vgamma}
(T_0+\frac{b c}{q})(T_0+1)=0,\\
\label{sahi4-Vgamma}
T_1 X - a =W(T_1+1-a^2),\\
\label{sahi6-Vgamma}
q T_0 W + c = X (T_0+1+\frac{b c}{q}),
\eea
which is the image of $\mathcal H_{V}$ under the isomorphism:
\bea\label{eq:gamma-V-sahi}
&&
\gamma(T_0,T_1,X,W)= (b T_1^{-1} W, T_1, i \sqrt{a b} T_1^{-1} W (T_0+1) X,-\frac{i}{\sqrt{a b}} W T_0 X T_1),\\
&&
\gamma(a,b,c) = (i\sqrt{a b}, i c\sqrt{\frac{b}{a}},- i \frac{q}{c}\sqrt{\frac{b}{a}}).\nn
\eea
Note that this fact has an interesting consequence in terms of $q$--polynomials: we shall see in Section \ref{se:zhedanov-others} that the spherical sub-algebra of  $\mathcal H_{V}$  acts as symmetries both on the continuous dual $q$--Hahn polynomials and the big $q$--Jacobi polynomials.

\subsection{Derivation of $\mathcal H_{IV}$,  $\mathcal H_{III}$, $\mathcal H_{II}$, $\mathcal H_{I}$}\label{der-1}

Again, the derivation procedure is more transparent in Noumi--Stokman representation. In each case we do a rescaling and  then take the limit $\eps\to 0$:

The algebra $\mathcal H_{IV}$ is obtained  from  $\mathcal H_{V}^\gamma$  by rescaling $V_0\to\frac{1}{\eps} V_0$ and $\check{V_0}\to\frac{1}{\eps} \check{V_0}$, $k_1\to\eps$, and $u_0\to\eps u_0$. 

The algebra $\mathcal H_{III}$ is obtained  from  $\mathcal H_{V}$ by rescaling $V_0\to\frac{1}{\eps} V_0$, $\check{V_0}\to\frac{1}{\eps}\check{V_0}$  and $u_0\to{\eps}\sqrt{q}$. 

The algebra $\mathcal H_{II}$ is obtained  from $\mathcal H_{IV}$ by rescaling $\check{V_0}\to\frac{1}{\eps}\check{V_0}$, $V_1\to\frac{1}{\eps}V_1$ $u_0\to\eps$. 

The algebra  $\mathcal H_{I}$ is obtained from $\mathcal H_{II}$ by rescaling $\check{V_1}\to\frac{1}{\eps} \check{V_1}$, $\check{V_0}\to\frac{1}{\eps}\check{V_0}$ and $u_1\to{\eps}$.

 \subsection{Derivation of  $\mathcal H_{III^{D_7}}$, and $\mathcal H_{III^{D_8}}$}\label{der-2}
 
In this case we don't have a Noumi--Stokman representation. We start form $\mathcal H_{III}$ in the presentation (\ref{sahi1-III}\dots,\ref{sahi6-III}). Rewrite relation (\ref{sahi4-III}) by using (\ref{sahi2-III}):
$$
T_1 X +(a+b) = W (T_1-(a+b)).
 $$
 Then take the limit as $b\to 0$ to obtain (\ref{sahiPIIID71},\dots\ref{sahiPIIID75}).
 
 Analogously, starting from (\ref{sahiPIIID71},\dots\ref{sahiPIIID75}) and taking $a\to 0$, we obtain the  $\mathcal H_{III^{D_8}}$ algebra.

\section{Confluent spherical sub--algebras and Painlev\'e cubics}\label{se:spherical-others}

In this section we give the confluent version of the results of section \ref{suse:cubic} for each algebra $\mathcal H_V,\mathcal H_{IV},\mathcal H_{III},
\mathcal H_{III^{D_7}},\mathcal H_{III^{D_8}},\mathcal H_{II},\mathcal H_{I}$ to produce in each case the appropriate cubic surface as limit of the spherical subalgebras. 
In each case we will select a idempotent element $e$ and define the spherical sub-algebra to be $e\mathcal H_{d} e$. We will prove all results in the first case, namely $\mathcal H_V$ and its spherical sub-algebra. The proofs in all other cases follow the same principles and we omit them for brevity.

Throughout this section we assume $q^m\neq 1$ for all $m\in\mathbb Z$.

\subsection{Spherical sub--algebra of $\mathcal H_V$ and PV monodromy manifold}

\begin{prop}\label{pr:cubic-V}
The following three elements:
\bea\label{sphericalPV}
\noindent& X_1 =X+W,\qquad 
X_2 =i \sqrt{a b}(T_0+1)T_1^{-1} -i \sqrt{\frac{1}{a b}} T_1T_0,\\
& X_3 =-i \sqrt{a b q} \left(\frac{1}{a b}W T_0 T_1 +\frac{1}{a b}X T_1 T_0+
\left(\frac{1}{a}+\frac{1}{b}\right) T_0+\frac{c}{a b q} T_1
\right),
\nn\eea
commute with $e=\frac{1+T_1}{1-a b}$,  satisfy the following  quantum commutation relations:
\bea\label{q-comm-PV}
&&q^{\frac{1}{2}} X_2 X_1-q^{-\frac{1}{2}} X_1 X_2 =  \left(q-\frac{1}{q}\right)X_3 +\nn\\
&&\qquad+  \left(q^{\frac{1}{2}} -q^{-\frac{1}{2}} \right)\left(\overline k_1+\frac{c}{\sqrt{q}}\left(\frac{i}{\sqrt{q a b}} T_1 + i \sqrt{q a b} T_1^{-1}\right)\right),\\
&&q^{\frac{1}{2}} X_3 X_2-q^{-\frac{1}{2}} X_2 X_3 =  \left(q^{\frac{1}{2}} -q^{-\frac{1}{2}} \right) \frac{c}{\sqrt{q}},\nn\\
&&q^{\frac{1}{2}} X_1 X_3-q^{-\frac{1}{2}} X_3 X_1 =  \left(q-\frac{1}{q}\right)X_2 -\nn\\
&&\qquad- \left(q^{\frac{1}{2}} -q^{-\frac{1}{2}} \right)\left(\overline k_1\frac{c}{\sqrt{q}}+\left(\frac{i}{\sqrt{q a b}} T_1 + i \sqrt{q a b} T_1^{-1}\right)\right),\nn
\eea
and the quantum cubic relation:
\bea\label{eq:skein-cubicV}
&&
q^{\frac{1}{2}}X_2 X_1 X_3-q X_2^2-q X_3^2+\sqrt{q} \left(-\overline k_1-\frac{c}{\sqrt{q}}\left(\frac{i}{\sqrt{q a b}} T_1 + i \sqrt{q a b} T_1^{-1}\right)\right) X_3-\nn\\
&&-  \frac{c}{q}X_1+ \sqrt{q}
\left(\overline k_1\frac{c}{\sqrt{q}}+\left(\frac{i}{\sqrt{q a b}} T_1 + i \sqrt{q a b} T_1^{-1}\right)\right)X_2+\nn\\
&&
\qquad \qquad +1+\frac{c^2}{q}- \frac{c}{\sqrt{q}}\overline k_1 \left(\frac{i}{\sqrt{q a b}} T_1 + i \sqrt{q a b} T_1^{-1}\right)=0,
\eea
where $ \overline k_1=i\sqrt{\frac{b}{a}}+i\sqrt{\frac{a}{b}}$.
\end{prop}

\proof
To prove the quantum commutation relations and the cubic relation we use the following substitutions which can be easily be derived from (\ref{sahi2-V},\ref{sahi3-V},\ref{sahi4-V}) and (\ref{sahi6-V}):
\begin{eqnarray}\label{q-alg-orderingV}
T_1^2=-(a b+1) T_1-a b, \qquad T_0^2=- T_0 \nn\\
T_1 X=W T_1+(1+ a b )W-(a+b),\\
 T_1 W= X T_1-(1+a b) W + a+b,\nn\\
T_0 X= q W T_0+c- X,\nn\\
T_0 W =\frac{1}{q} X T_0+\frac{1}{q}X-\frac{c}{q}.\nn
\end{eqnarray}
Using these relations recursively, it is a straightforward computation to eliminate all higher powers in $T_0$ and $T_1$ and to bring all $X$ and $W$ to the left and eventually obtain the statement.
\endproof

\begin{cor}\label{co:cubic-V}
In the  limit $q\to 1$ the elements $X_1,X_2,X_3$ belong to the centre of $\mathcal H_{V}$ and the cubic relation (\ref{eq:skein-cubicV}) tends to  the PV monodromy manifold  (see Appendix A):
\be\label{cubic-PV}
X_1 X_2 X_3 -X_2^2-X_3^2-c  X_1- \left( \overline u_1-  {c \overline k_1}\right) X_2-\left(-c \overline u_1 +\overline k_1\right) X_3 +1+ {c^2}+ c \overline k_1   \overline u_1=0.
\ee
\end{cor}

\proof
Observe that in the  limit $q\to 1$,
$$
\frac{i}{\sqrt{q a b}} T_1 + i \sqrt{q a b} T_1^{-1}\to -\frac{i(1+a b)}{\sqrt{ a b}} = -i\overline u_1,
$$
so that we obtain (\ref{cubic-PV}) by substituting $q=1$ and $\frac{i}{\sqrt{q a b}} T_1 + i \sqrt{q a b} T_1^{-1}\to -\frac{i(1+a b)}{\sqrt{ a b}}$ by $  -i\overline u_1$  in (\ref{eq:skein-cubicV}).
\endproof

\begin{cor}\label{pr:spherical-PV}
Let $e=\frac{1+T_1}{1-a b}$, then
the elements $\hat X_i= e X_i e$, $i=1,2,3$, where $X_1,X_2,X_3$ are defined by (\ref{sphericalPV}), generate the spherical sub-algebra $e\mathcal H_{V}e$, they satisfy the quantum commutation relations:
\bea\label{eq:skein-hat-5}
&&q^{\frac{1}{2}}\hat X_2 \hat X_1-q^{-\frac{1}{2}} \hat X_1 \hat X_2 =  \left(q-\frac{1}{q}\right)\hat X_3 - 
 \left(q^{\frac{1}{2}} -q^{-\frac{1}{2}} \right)\omega_3 e,\nn\\
&&q^{\frac{1}{2}} \hat X_3 \hat X_2-q^{-\frac{1}{2}} \hat X_2 \hat X_3 =  - 
 \left(q^{\frac{1}{2}} -q^{-\frac{1}{2}} \right)\omega_1 e,\\
&&q^{\frac{1}{2}} \hat X_1 \hat X_3-q^{-\frac{1}{2}} \hat X_3 \hat X_1 =  \left(q-\frac{1}{q}\right)\hat X_2 -
 \left(q^{\frac{1}{2}} -q^{-\frac{1}{2}} \right)\omega_2 e,\nn
 \eea
and the quantum cubic relation:
\be\label{eq:cubic-hatV}
q^{\frac{1}{2}} \hat X_2 \hat X_1\hat X_3 -
q \hat X_2^2- q \hat X_3^2+\sqrt{q} \omega_2 \hat X_2+ \frac{1}{\sqrt{q}} \omega_1 \hat X_1+\sqrt{q} \omega_3 \hat X_3-\omega_4 e=0,
\ee
where
\begin{eqnarray}\label{omegaPV}
&&
\omega_1=- \frac{c}{\sqrt{q}},\qquad \omega_2= \frac{c}{\sqrt{q}}\overline k_1 -i \left(\sqrt{\frac{a b}{q}}+\sqrt{\frac{q}{a b}}\right),\nn\\
&&
\omega_3=-\overline k_1+ i \frac{c}{\sqrt{q}}
 \left(\sqrt{\frac{a b}{q}}+\sqrt{\frac{q}{a b}}\right),\\
&&
\omega_4=1+\frac{c^2}{q}+i \frac{c \overline k_1}{ \sqrt{q}}  \left(\sqrt{\frac{a b}{q}}+\sqrt{\frac{q}{a b}}\right).\nn
\end{eqnarray}
\end{cor}

\proof
The fact that $\hat X_1,\hat X_2,\hat X_3$ and $e$ generate the spherical sub-algebra $e\mathcal H_{V}$  follows easily from the following relations
$$
eT_0e=\frac{i\sqrt{a b}}{1-ab}\hat X_2-\frac{1}{1-a b} e, \quad e T_1 e=- a b e,\quad 
e X e= \frac{1}{1-a b}\hat X_1-\frac{a+b}{1-a b}e,
$$
which can be proved as in the proof of Proposition \ref{pr:cubic-V}.

To prove the quantum commutation relations (\ref{eq:skein-hat-5}) it is enough to observe that $e$ is idempotent and to prove that $X_1,X_2,X_3$ commute with $e$. Indeed if this is true, we can just multiply (\ref{q-comm-PV}) by $e$ and use the fact that $e X_i X_j = e^2 X_i X_j= e X_i e X_j= \hat X_i \hat X_j$. In a similar way we can prove (\ref{eq:cubic-hatV}) by multiplying (\ref{eq:skein-cubicV}) by $e$ three times.

So we only need to prove that $[e,X_{1,2}]=0$:
\bea
&&
(1-a b)[e,X_1] = [1+T_1, X+W] =T_1 X +T_1 W -X T_1 - W T_1= \nn\\
&&
=(W T_1+(1+ a b )W-(a+b))+(X T_1-(1+a b) W + a+b)-X T_1 - W T_1=0.\nn\\
&&
(1-a b)[e,X_2] = [1+T_1,i \sqrt{a b}(T_0+1)T_1^{-1} -i \sqrt{\frac{1}{a b}} T_1T_0] = \nn\\
&&
= i \sqrt{a b}T_1 (T_0+1)T_1^{-1} -i \sqrt{\frac{1}{a b}} T_1^2T_0- i \sqrt{a b}T_0(T_0+1)-i \sqrt{\frac{1}{a b}} T_0 T_1T_0=\nn\\
&&
=- i \sqrt{a b}T_1 (T_0+1)(\frac{1}{ab} T_1+\frac{1+a b}{a b}) -i \sqrt{\frac{1}{a b}} (-(a b+1) T_1-a b) T_0 -i \sqrt{\frac{1}{a b}} T_0 T_1T_0=0.
\nn\eea
This concludes the proof.
\endproof

\subsection{Spherical sub--algebra of $\mathcal H_{IV}$ and PIV monodromy manifold}

\begin{prop}\label{pr:cubic-IV} 
The following three elements:
\bea\label{sphericalPIV}
\noindent& X_1 =X+W,\qquad 
X_2 =i \sqrt{a b}(T_0+1)T_1^{-1} -i \sqrt{\frac{1}{a b}} T_1T_0,\\
& X_3 =-i \sqrt{\frac{q}{a b} } \left( W T_0 T_1 + X T_1 T_0-i\sqrt{a b}
T_0+\frac{c}{ q} T_1\right)
\eea
commute with $e=\frac{1+T_1}{1-a b}$,  satisfy the following  quantum commutation relations:
\bea
&&q^{\frac{1}{2}} X_2 X_1-q^{-\frac{1}{2}} X_1 X_2 =  \left(q-\frac{1}{q}\right)X_3 +\nn\\
&&\qquad+  \left(q^{\frac{1}{2}} -q^{-\frac{1}{2}} \right)
\left(-1+ \frac{c}{\sqrt{q}}\left(\frac{i}{\sqrt{q a b}} T_1 + i \sqrt{q a b} T_1^{-1}\right)\right),\nn\\
&&q^{\frac{1}{2}} X_3 X_2-q^{-\frac{1}{2}} X_2 X_3 =   \left(q^{\frac{1}{2}} -q^{-\frac{1}{2}} \right)  \frac{c}{\sqrt{q}},\nn\\
&&q^{\frac{1}{2}} X_1 X_3-q^{-\frac{1}{2}} X_3 X_1 =  \left(q^{\frac{1}{2}} -q^{-\frac{1}{2}} \right)\frac{c}{\sqrt{q}},\nn
\eea
and the quantum cubic relation:
\bea\label{eq:skein-cubicIV}
&&
q^{\frac{1}{2}}X_2 X_1 X_3-q X_3^2-c X_2- \frac{c}{q}X_1
+\frac{c^2}{q}+\frac{c}{\sqrt{q}}\left(\frac{i}{\sqrt{q a b}} T_1 + i \sqrt{q a b} T_1^{-1}\right)+\nn\\
&&\qquad\qquad + \sqrt{q}\left(1- \frac{c}{\sqrt{q}}\left(\frac{i}{\sqrt{q a b}} T_1 + i \sqrt{q a b} T_1^{-1}\right)\right) X_3.
\eea
\end{prop}

\begin{cor}\label{co:cubic-IV}
In the limit $q\to 1$ the elements $X_1,X_2,X_3$ belong to the centre of $\mathcal H_{IV}$ and the cubic relation (\ref{eq:skein-cubicIV}) tends to the PIV monodromy manifold  (see Appendix A):
\be\label{cubic-PIV}
X_1 X_2 X_3-X_3^2- c   X_1-c  X_2+( 1+\frac{c}{b}-c b) X_3 + c(c+b-\frac{1}{b})=0.
\ee
\end{cor}

\begin{cor}\label{pr:spherical-PIV}
The elements $\hat X_i= e X_i e$, $i=1,2,3$, where $X_1,X_2,X_3$ are defined by (\ref{sphericalPIV}), generate the spherical sub-algebra $e\mathcal He$, satisfy the quantum commutation relations:
\bea\label{eq:skein-hat-4}
&&q^{\frac{1}{2}}\hat X_2 \hat X_1-q^{-\frac{1}{2}} \hat X_1 \hat X_2 =  \left(q-\frac{1}{q}\right)\hat X_3 - 
 \left(q^{\frac{1}{2}} -q^{-\frac{1}{2}} \right)\omega_3 e,\nn\\
&&q^{\frac{1}{2}} \hat X_3 \hat X_2-q^{-\frac{1}{2}} \hat X_2 \hat X_3 =  - 
 \left(q^{\frac{1}{2}} -q^{-\frac{1}{2}} \right)\omega_1 e,\\
&&q^{\frac{1}{2}} \hat X_1 \hat X_3-q^{-\frac{1}{2}} \hat X_3 \hat X_1 =  -
 \left(q^{\frac{1}{2}} -q^{-\frac{1}{2}} \right)\omega_2 e,\nn
 \eea
 and the quantum cubic relation
 \bea\label{eq:cubic-hatIV}
&&
q^{\frac{1}{2}} \hat X_2 \hat X_1\hat X_3 - q \hat X_3^2+\sqrt{q} \omega_2 \hat X_2+ \frac{1}{\sqrt{q}} \omega_1 \hat X_1+\sqrt{q} \omega_3 \hat X_3+
\omega_4 e=0,
\eea
where
\be\label{omegaPIV}
\omega_1=\omega_2=-\frac{c}{\sqrt{q}},\quad
\omega_3= 1+\frac{c}{\sqrt{q}}  i \left(\sqrt{\frac{a b}{q}}+\sqrt{\frac{q}{a b}}\right),\quad
\omega_4=\frac{c^2}{q}-\frac{c}{\sqrt{q}}  i \left(\sqrt{\frac{a b}{q}}+\sqrt{\frac{q}{a b}}\right).
\ee
\end{cor}

\subsection{Spherical sub--algebra of $\mathcal H_{III}$ and PIII monodromy manifold}

\begin{prop}\label{pr:cubic-III}
The following three elements:
\bea\label{sphericalPIII}
\noindent& X_1 =X+W,\qquad 
X_2 =i \sqrt{a b}T_0 T_1^{-1} -i \sqrt{\frac{1}{a b}} T_1T_0,\\
& X_3 =-i \sqrt{a b q} \left(\frac{1}{a b}W T_0 T_1 +\frac{1}{a b}X T_1 T_0+
\left(\frac{1}{a}+\frac{1}{b}\right) T_0-\frac{1}{a b q} T_1
\right),
\eea
commute with $e=\frac{1+T_1}{1-a b}$,  satisfy the following  quantum commutation relations:
\bea
&&q^{\frac{1}{2}} X_2 X_1-q^{-\frac{1}{2}} X_1 X_2 =  \left(q-\frac{1}{q}\right)X_3 -
q^{-\frac{1}{2}}   \left(q^{\frac{1}{2}} -q^{-\frac{1}{2}} \right) \left(\frac{i}{\sqrt{q a b}} T_1 + i \sqrt{q a b} T_1^{-1}\right),\nn\\
&&q^{\frac{1}{2}} X_3 X_2-q^{-\frac{1}{2}} X_2 X_3 =0,\nn\\
&&q^{\frac{1}{2}} X_1 X_3-q^{-\frac{1}{2}} X_3 X_1 =  \left(q-\frac{1}{q}\right)X_2 +\left(q^{\frac{1}{2}} -q^{-\frac{1}{2}} \right)q^{-\frac{1}{2}} \overline k_1,\nn
\eea
and the quantum cubic relation:
\be\label{eq:skein-cubicIII}
q^{\frac{1}{2}}X_2 X_1 X_3-q X_2^2-q X_3^2- \overline k_1X_2 + \left(\frac{i}{\sqrt{q a b}} T_1 + i \sqrt{q a b} T_1^{-1}\right) X_3+1/q=0,
\ee
where $ \overline k_1= i\sqrt{\frac{b}{a}}+i\sqrt{\frac{a}{b}}$.
\end{prop}

\begin{cor}\label{co:cubic-III}
In the limit $q\to 1$ the elements $X_1,X_2,X_3$ generate the centre of $\mathcal H_{III}$ and the cubic relation (\ref{eq:skein-cubicIII}) tends to the PIII monodromy manifold  (see Appendix A):
\be\label{cubic-PIII}
X_1 X_2 X_3 -X_2^2-X_3^2-\left( i\sqrt{\frac{b}{a}}+i\sqrt{\frac{a}{b}}\right) X_2-\left(i\sqrt{a b} +i\frac{1}{\sqrt{a b}}\right)  X_3 +1=0.
\ee
\end{cor}

\begin{cor}\label{pr:spherical-PIII}
The elements $\hat X_i= e X_i e$, $i=1,2,3$, where $X_1,X_2,X_3$ are defined by (\ref{sphericalPIII}), generate the spherical sub-algebra $e\mathcal He$, they satisfy the quantum commutation relations:
\bea\label{eq:skein-hat-3}
&&q^{\frac{1}{2}}\hat X_2 \hat X_1-q^{-\frac{1}{2}} \hat X_1 \hat X_2 =  \left(q-\frac{1}{q}\right)\hat X_3 - 
 \left(q^{\frac{1}{2}} -q^{-\frac{1}{2}} \right)\omega_3 e,\nn\\
&&q^{\frac{1}{2}} \hat X_3 \hat X_2-q^{-\frac{1}{2}} \hat X_2 \hat X_3 =  0,\\
&&q^{\frac{1}{2}} \hat X_1 \hat X_3-q^{-\frac{1}{2}} \hat X_3 \hat X_1 =  \left(q-\frac{1}{q}\right)\hat X_2 -
 \left(q^{\frac{1}{2}} -q^{-\frac{1}{2}} \right)\omega_2 e,\nn
 \eea
 and the quantum cubic relation:
\be\label{eq:cubic-hatIII}
q^{\frac{1}{2}} \hat X_2 \hat X_1\hat X_3 -
q \hat X_2^2- q \hat X_3^2+\sqrt{q} \omega_2 \hat X_2+\sqrt{q} \omega_3 \hat X_3+ e=0,
\ee 
where
\begin{eqnarray}\label{omegaPIII}
&&
\omega_2=-\frac{i}{\sqrt{q}}\left(\sqrt{\frac{b}{a}}-i\sqrt{\frac{a}{b}}\right),\qquad \omega_3= -\frac{i}{\sqrt{q}} \left(\sqrt{\frac{a b}{q}}+\sqrt{\frac{q}{a b}}\right).\nn
\end{eqnarray}
\end{cor}

\subsection{Spherical sub--algebra of $\mathcal H_{III^{D_7}}$ and $PIII^{D_7}$ monodromy manifold}

\begin{lm}
The generators $X, W,T_0,T_1$ of the confluent Cherednik algebra $\mathcal H_{III}^{D_7}$ satisfy the following relations
\bea
\label{sahiPIIID74-1}
&&
T_1 X T_1  = - a T_1,\qquad
(T_1+1) W (T_1+1)  = a (T_1+1),\\
\label{sahiPIIID75-1}
&&
T_0 W T_0 =-\frac{1}{q} T_0,\qquad T_0 X T_0=T_0.
\eea
\end{lm}

\begin{prop} \label{pr:cubic-IIID7}
The following three elements:
\bea\label{sphericalPIIID7}
X_1 = X+W,\nn\\
X_2 =T_1 T_0 +T_0(T_1+1)\\
X_3 =\frac{q}{q^2-1}\left(q^{1/2}X_2 X_1 -q^{-1/2} X_1 X_2\right) - \frac{1}{q+1}\left(\left(\sqrt{q}-\frac{1}{\sqrt{q}}\right)T_1+\sqrt{q}\right),\nn
\eea
commute with $e=1+T_1$, they satisfy the quantum commutation relations:
\bea\label{skein-PIIID7}
&&
q^{\frac{1}{2}} X_2 X_1-q^{-\frac{1}{2}} X_1 X_2 =  \left(q-\frac{1}{q}\right)X_3 + \frac{q-1}{q}\left(\left(\sqrt{q}-\frac{1}{\sqrt{q}}\right)T_1+\sqrt{q}\right),\nn\\
&&
q^{\frac{1}{2}} X_3 X_2-q^{-\frac{1}{2}} X_2 X_3 =  0,\\
&&
q^{\frac{1}{2}} X_1 X_3-q^{-\frac{1}{2}} X_3 X_1 = (q-\frac{1}{q})X_2 -\frac{q-1}{q}a,\nn
\eea
and the quantum cubic relation:
\be\label{eq:skein-cubicIIID7}
q^{\frac{1}{2}}X_2 X_1 X_3-q X_2^2-q X_3^2+ a X_2 +(q^{-\frac{1}{2}} T_1 -q^{\frac{1}{2}} (T_1+1)) X_3
=0.
\ee
\end{prop}

\begin{cor}\label{co:cubic-IIID7}
In the limit $q\to 1$, the elements $X_1,X_2,X_3$ become central and the cubic relation  (\ref{eq:skein-cubicIIID7}) tends to the the $P{III}^{D_7}$ monodromy manifold  (see Appendix A):
\be\label{eq:D7cubicq1}
X_3 X_2 X_1-X_2^2-  X_3^2 -a\, X_2+ X_3=0.
\ee
\end{cor}

\begin{cor}\label{pr:spherical-PIIID7}
Define $\hat X_i= e X_i e$, $i=1,2,3$, where $X_1,X_2,X_3$ are defined by (\ref{sphericalPIIID7}) and 
$$
e=1+T_1.
$$
Then $\hat X_1$, $\hat X_2$,  $\hat X_3$ generate the spherical sub-algebra $e\mathcal H_{III}^{D_7}e$, they satisfy the quantum commutation relations:
\bea\label{eq:D7-skein}
&&
q^{\frac{1}{2}}\hat X_2 \hat X_1-q^{-\frac{1}{2}} \hat X_1 \hat X_2 =  \left(q-\frac{1}{q}\right) \hat X_3 - \left(q^{\frac{1}{2}} -q^{-\frac{1}{2}} \right)e,  \nn \\
&&
q^{\frac{1}{2}} \hat X_3 \hat X_2-q^{-\frac{1}{2}} \hat X_2 \hat X_3 =  0,\\
&&
q^{\frac{1}{2}} \hat X_1 \hat X_3-q^{-\frac{1}{2}} \hat X_3 \hat X_1 =  \left(q-\frac{1}{q}\right)\hat X_2  
- \left(q^{\frac{1}{2}} -q^{-\frac{1}{2}} \right)\frac{a}{\sqrt{q}}\, e, \nn
\eea
and the quantum cubic cubic relation:
\be\label{eq:D7cubic}
q^{\frac{1}{2}} \hat X_2 \hat X_1\hat X_3 -
q \hat X_2^2- q \hat X_3^2+ a \hat X_2-\frac{1}{\sqrt{q}}   \hat X_3=0.
\ee
\end{cor}

\subsection{Spherical sub--algebra of $\mathcal H_{III^{D_8}}$ and $PIII^{D_8}$ monodromy manifold}

Here all proofs are a simple limit as $a\to 0$ of the proofs of the previous Sub-section and will be omitted.
\begin{lm}
The generators $X,W,T_0,T_1$ of the confluent Cherednik algebra $\mathcal H_{III^{D_8}}$ satisfy the following relations
\bea
\label{sahiPIIID84-1}
&&
T_1 X T_1  =0,\qquad
(T_1+1) W (T_1+1)  = 0,\\
\label{sahiPIIID85-1}
&&
T_0 W T_0 =-\frac{1}{q} T_0,\qquad T_0 X T_0=T_0.
\eea
\end{lm}

\begin{prop} \label{pr:cubic-IIID8}
The following three elements:
\bea\label{sphericalPIIID8}
X_1 = X+W,\nn\\
X_2 =T_1 T_0 +T_0(T_1+1)\\
X_3 =\frac{q}{q^2-1}\left(q^{1/2}X_2 X_1 -q^{-1/2} X_1 X_2\right) - \frac{1}{q+1}\left(\left(\sqrt{q}-\frac{1}{\sqrt{q}}\right)T_1+\sqrt{q}\right),\nn
\eea
commute with $e=1+T_1$,  satisfy the following  quantum commutation relations:
\bea\label{skein-PIIID8}
&&
q^{\frac{1}{2}} X_2 X_1-q^{-\frac{1}{2}} X_1 X_2 =  \left(q-\frac{1}{q}\right)X_3 + \frac{q-1}{q}\left(\left(\sqrt{q}-\frac{1}{\sqrt{q}}\right)T_1+\sqrt{q}\right),\nn\\
&&
q^{\frac{1}{2}} X_3 X_2-q^{-\frac{1}{2}} X_2 X_3 =  0,\\
&&
q^{\frac{1}{2}} X_1 X_3-q^{-\frac{1}{2}} X_3 X_1 = (q-\frac{1}{q})X_2,\nn
\eea
and the quantum cubic relation:
\be\label{eq:skein-cubicIIID8}
q^{\frac{1}{2}}X_2 X_1 X_3-q X_2^2-q X_3^2 + (q^{-\frac{1}{2}}T_1 -q^{\frac{1}{2}} (T_1+1)) X_3
=0.
\ee
\end{prop}

\begin{cor}\label{co:cubic-IIID8}
In the limit $q\to 1$, the elements $X_1,X_2,X_3$ become central and the cubic relation (\ref{eq:skein-cubicIIID8}) tends to the the $P{III}^{D_8}$ monodromy manifold  (see Appendix A):
\be\label{eq:D8cubicq1}
X_3 X_2 X_1-X_2^2-  X_3^2+ X_3=0.
\ee
\end{cor}

\begin{cor}\label{pr:spherical-PIIID8}
Define $\hat X_i= e X_i e$, $i=1,2,$, where $X_1,X_2,X_3$ are defined by (\ref{sphericalPIIID8}) and 
$$
e=1+T_1.
$$
Then $\hat X_1$, $\hat X_2$ generate the spherical sub-algebra $e\mathcal H_{III^{D_8}}e$, the satisfy the quantum commutation relations:
\bea\label{eq:D8-skein}
&&
q^{\frac{1}{2}}\hat X_2 \hat X_1-q^{-\frac{1}{2}} \hat X_1 \hat X_2 =  \left(q-\frac{1}{q}\right) \hat X_3 - \frac{q-1}{\sqrt{q}} e,  \nn \\
&&
q^{\frac{1}{2}} \hat X_3 \hat X_2-q^{-\frac{1}{2}} \hat X_2 \hat X_3 =  0,\\
&&
q^{\frac{1}{2}} \hat X_1 \hat X_3-q^{-\frac{1}{2}} \hat X_3 \hat X_1 =  \left(q-\frac{1}{q}\right)\hat X_2, \nn
\eea
and lie on the following quantum cubic
\be\label{eq:D8cubic}
q^{\frac{1}{2}} \hat X_2 \hat X_1\hat X_3 -
q \hat X_2^2- q \hat X_3^2-\frac{1}{\sqrt{q}}   \hat X_3=0.
\ee
\end{cor}

\subsection{Spherical sub--algebra of $\mathcal H_{II}$ and PII monodromy manifold}

\begin{prop}\label{pr:cubic-II}
The following three elements:
\bea\label{sphericalPII}
\noindent& X_1 =X+W,\qquad 
X_2 = - b (T_0+1) T_1^{-1} -\frac{1}{ b} T_1T_0,\\
& X_3 =-\frac{\sqrt{q}}{b} \left( W T_0 T_1 + X T_1 T_0-\frac{1}{\sqrt{ q}} T_1\right)
\eea
commute with $e=\frac{1+T_1}{1+b^2}$,  satisfy the following  quantum commutation relations:
\bea
&&q^{\frac{1}{2}} X_2 X_1-q^{-\frac{1}{2}} X_1 X_2 =  \left(q-\frac{1}{q}\right)  X_3
 +\left(q^{\frac{1}{2}} -q^{-\frac{1}{2}} \right)
\left(-\frac{1}{\sqrt{q} b} T_1 + \sqrt{q} b T_1^{-1}\right),\nn\\
&&q^{\frac{1}{2}} X_3 X_2-q^{-\frac{1}{2}} X_2 X_3 =  -\left(q^{\frac{1}{2}} -q^{-\frac{1}{2}} \right),\nn\\
&&q^{\frac{1}{2}} X_1 X_3-q^{-\frac{1}{2}} X_3 X_1 = 0,\nn
\eea
and the quantum cubic relation:
\be\label{eq:skein-cubicII}
q^{\frac{1}{2}}X_2 X_1 X_3-q X_3^2 + \frac{1}{\sqrt{q}}  X_1+
 \sqrt{q} \left(-\frac{1}{\sqrt{q} b} T_1 + \sqrt{q} b T_1^{-1}\right)X_3+1=0
\ee\end{prop}

\begin{cor}\label{co:cubic-II}
In the limit $q\to 1$ the elements $X_1,X_2,X_3$ belong to the centre of $\mathcal H_{II}$ and the cubic relation (\ref{eq:skein-cubicII}) tends to the PII monodromy manifold  (see Appendix A):
\be\label{cubic-PII}
X_1 X_2 X_3 -X_3^2+ X_1+\left(b-\frac{1}{b}\right)X_3 +1 =0.
\ee
\end{cor}

\begin{cor}\label{pr:spherical-PII}
The elements $\hat X_i= e X_i e$, $i=1,2,3$, where $X_1,X_2,X_3$ are defined by (\ref{sphericalPII}), generate the spherical sub-algebra $e\mathcal H_{II} e$, satisfy the quantum commutation relations:
\bea\label{eq:skein-hat-2}
&&q^{\frac{1}{2}}\hat X_2 \hat X_1-q^{-\frac{1}{2}} \hat X_1 \hat X_2 =   \left(q-\frac{1}{q}\right) \hat X_3 - 
 \left(q^{\frac{1}{2}} -q^{-\frac{1}{2}} \right)\omega_3 e,\nn\\
&&q^{\frac{1}{2}} \hat X_3 \hat X_2-q^{-\frac{1}{2}} \hat X_2 \hat X_3 =  -
 \left(q^{\frac{1}{2}} -q^{-\frac{1}{2}} \right) e,\\
&&q^{\frac{1}{2}} \hat X_1 \hat X_3-q^{-\frac{1}{2}} \hat X_3 \hat X_1 = 0,\nn
 \eea
 and the quantum cubic relation
 \bea\label{eq:cubic-hatII}
&&
q^{\frac{1}{2}} \hat X_2 \hat X_1\hat X_3-q \hat X_3^2 +\frac{1}{\sqrt{q}} \hat X_1+\sqrt{q} \omega_3 \hat X_3+ e=0,
\eea
where
\be\label{omegaPII}
\omega_3=\left(\frac{b}{\sqrt{q}}-\frac{\sqrt{q}}{b}\right).
\ee
\end{cor}

\subsection{Spherical sub--algebra of  $\mathcal H_{I}$ and PI monodromy manifold}

\begin{prop}\label{pr:cubic-I}
The following three elements:
\bea\label{sphericalPI}
X_1 = \check{V_1} V_1+V_1(\check{V_1}+1),\nn\\
X_2 =\check{V_1} V_0 + (V_0+1) (\check{V_1}+1)\\
X_3 =q^{1/2} {V_1} V_0+q^{-1/2} (V_0+1)V_1,\nn
\eea
commute with $e=1+\check V_1$,  satisfy the following  quantum commutation relations:
\bea
&&q^{\frac{1}{2}} X_2 X_1-q^{-\frac{1}{2}} X_1 X_2 = 
  \left(q^{\frac{1}{2}} -q^{-\frac{1}{2}} \right)
(q^{-\frac{1}{2}} \check{ V_1} -q^{\frac{1}{2}}( \check{ V_1}+1)),\nn\\
&&q^{\frac{1}{2}} X_3 X_2-q^{-\frac{1}{2}} X_2 X_3 = -\left(q^{\frac{1}{2}} -q^{-\frac{1}{2}} \right),\nn\\
&&q^{\frac{1}{2}} X_1 X_3-q^{-\frac{1}{2}} X_3 X_1 =  0,\nn
\eea
and the quantum cubic relation:
\be\label{eq:skein-cubicI}
q^{\frac{1}{2}}X_2 X_1 X_3+ \frac{1}{\sqrt{q}}   X_1- \sqrt{q} 
(q^{-\frac{1}{2}} \check{ V_1} -q^{\frac{1}{2}}( \check{ V_1}+1))X_3+1=0.
\ee
\end{prop}

\begin{cor}\label{co:cubic-I}
In the limit $q\to 1$ the elements $X_1,X_2,X_3$ belong to the centre of $\mathcal H_{I}$ and the cubic relation (\ref{eq:skein-cubicI}) tends to the PI monodromy manifold  (see Appendix A):
\be\label{cubic-PI}
X_1 X_2 X_3+   X_1+  X_3 +1=0.
\ee
\end{cor}

\begin{cor}\label{pr:spherical-PI}
The elements $\hat X_i= e X_i e$, $i=1,2,3$, where $X_1,X_2,X_3$ are defined by (\ref{sphericalPI}), generate the spherical sub-algebra $e\mathcal H_{I} e$, satisfy the quantum commutation relations:
\bea\label{eq:skein-hat-1}
&&q^{\frac{1}{2}}\hat X_2 \hat X_1-q^{-\frac{1}{2}} \hat X_1 \hat X_2 =  - 
 \left(q^{\frac{1}{2}} -q^{-\frac{1}{2}} \right)\sqrt{q} \,e,\nn\\
&&q^{\frac{1}{2}} \hat X_3 \hat X_2-q^{-\frac{1}{2}} \hat X_2 \hat X_3 = - \left(q^{\frac{1}{2}} -q^{-\frac{1}{2}} \right) e,\\
&&q^{\frac{1}{2}} \hat X_1 \hat X_3-q^{-\frac{1}{2}} \hat X_3 \hat X_1 =  0,\nn
 \eea
 and the quantum cubic relation
 \bea\label{eq:cubic-hatI}
&&
q^{\frac{1}{2}} \hat X_2 \hat X_1\hat X_3 +\frac{1}{\sqrt{q}} \hat X_1+ q  \hat X_3+ e=0.
\eea
\end{cor}

\section{Embedding of the (confluent) Cherednik algebras into $Mat(2,\mathbb T_q)$}\label{se:embedding}

In this section we give an embedding for  $\mathcal H$, $\mathcal H_V$, $\mathcal H_{IV}$, $\mathcal H_{III}$ and 
$\mathcal H_{II}$  into  $Mat(2,\mathbb T_q)$. The proof of such embedding is based on the Lusztig-Demazure presentation given in the following:

\begin{lm}
The algebras $\mathcal H,\mathcal H_V,\mathcal H_{IV},\mathcal H_{III},\mathcal H_{II}$ are the 
 algebras generated by five elements $T^{\pm 1},X,W,Y,Z$, where
 \be\label{eq:LDty}
 T=T_1,\qquad Y=T_1 T_0,
\ee
 and
 \be\label{eq:LDz}
 Z= \left\{\begin{array}{lc}
(T_0+\frac{c d}{q}+1)(T_1+ a b+1),&\hbox{for}\quad \mathcal H,\\
(T_0+1)(T_1+ a b+1),&\hbox{for}\quad \mathcal H_{V}, \\
(T_0+1)(T_1- b^2+1),&\hbox{for}\quad \mathcal H_{IV},\mathcal H_{II},\\
T_0(T_1+ a b+1),&\hbox{for}\quad \mathcal H_{III},\\
\end{array}
 \right.
\ee
 satisfying the following relations respectively:
\begin{itemize}
\item $\mathcal H$
\bea
\label{LD0} X W =W X=1,\\
\label{LD00} Y Z =Z Y=1,\\
\label{LD1} X T + a b T^{-1} W +a+b=0,\\
\label{LD2} Z T +\frac{q}{c d} T^{-1} Y + 1+\frac{q}{c d}=0,\\
\label{LD3} (T+ a b)(T+1)=0,\\
\label{LD4} \qquad \quad Y X  =- \frac{q}{a b} T^2 X Y  -q\left( \frac{1}{a}+\frac{1}{b}\right) T Y - \left(1+\frac{c d}{q}\right)T X + (c+d)T.
\eea
\item $\mathcal H_{V}$
\bea
\label{LD0-PV} W X=X W=1,\\
\label{LD00-PV} Z Y=Y Z=0,\\
\label{LD1-PV} X T + a b T^{-1} W +a+b=0,\\
\label{LD2-PV} Z T + T^{-1} Y +1=0,\\
\label{LD3-PV}(T+ a b)(T+1)=0,\\
\label{LD4-PV} Y X  =- \frac{q}{a b} T^2 X Y  -q\left( \frac{1}{a}+\frac{1}{b}\right) T Y- T X + c T.
\eea
\item $\mathcal H_{IV}$:
\bea
\label{LD0-PIV} W X=X W=0,\\
\label{LD00-PIV} Z Y=Y Z=0,\\
\label{LD1-piv} X T = b^2 T^{-1} W-b,\\
\label{LD2-piv} Z T + T^{-1} Y+1=0,\\
\label{LD3-piv}  (T-b^2)(T+1)=0,\\
\label{LD4-piv} Y X  =\frac{q}{b^2} T^2 X Y  +\frac{q}{b} T Y -  T X + c T.
\eea
\item $\mathcal H_{III}$:
\bea
\label{LD0-PIII} W X=X W=1,\\
\label{LD00-PIII} Z Y=Y Z=0,\\
\label{LD1-piii} X T + a b T^{-1} W +a+b=0,\\
\label{LD2-piii} Z T +T^{-1} Y=0,\\
\label{LD3-piii} (T+ a b)(T+1)=0,\\
\label{LD4-piii} Y X  =- \frac{q}{a b} T^2 X Y  -q\left( \frac{1}{a}+\frac{1}{b}\right) T Y -  T.
\eea
\item $\mathcal H_{II}$:
\bea
\label{LD0-PII} W X=X W=0,\\
\label{LD00-PII} Z Y=Y Z=0,\\
\label{LD1-pii} X T = b^2 T^{-1} W,\\
\label{LD2-pii} Z T +T^{-1} Y+1=0,\\
\label{LD3-pii}  (T-b^2)(T+1)=0,\\
\label{LD4-pii} Y X  =\frac{q}{b^2} T^2 X Y   -  T X -\sqrt{q} T.
\eea
\end{itemize}
\end{lm}

\proof For $\mathcal H$, see  Proposition 6.6 in \cite{O}. For the other cases,  it is enough to observe that since $T$ is invertible, we can always invert relations (\ref{eq:LDty}) to define $T_0= T^{-1}Y$. To check that the algebra relations in the $X,W,T_0,T_1$ presentation imply the ones in the $X,W,T,Y,Z$ one and vice-versa is a straightforward computation.
\endproof

\begin{theorem}\label{th-rep-qT}
The map:
\be
\label{rep1}
T\to\left(\begin{array}{cc}
0&  \sqrt{a b}\, e^{S_1}\\
- \sqrt{a b}\, e^{-S_1}&-1-a b  \\
\end{array}\right),
\ee
\be
\label{rep2}
X\to\left(\begin{array}{cc}
\left(\sqrt{\frac{a}{b}}-\sqrt{\frac{b}{a}}\right)e^{-S_1}+q\, e^{-S_1}e^{-S_2} + \frac{1}{q} e^{-S_1}e^{S_2}-
\left(\sqrt{\frac{1}{a b}}+\sqrt{a b}\right) e^{S_2}& -q \, e^{S_1}e^{S_2} \\
-\left(\sqrt{\frac{a}{b}}-\sqrt{\frac{b}{a}}\right)e^{-S_1}- \frac{1}{q} e^{-S_1}e^{S_2}+
\left(\sqrt{\frac{1}{a b}}+\sqrt{a b}\right) e^{S_2}& q\, e^{S_1}e^{S_2}\\
\end{array}\right),
\ee

\bea
\label{rep3}
W&\to&\left(\begin{array}{cc}
e^{S_1}e^{S_2}& \\
\left(\sqrt{\frac{a}{b}}-\sqrt{\frac{b}{a}}\right)e^{-S_1}+e^{-S_1}e^{S_2} -\left(\sqrt{\frac{1}{a b}}+\sqrt{a b}\right) e^{S_2}& \\
\end{array}\right.\nn \\
&&\left.\qquad\qquad \qquad
\begin{array}{cc} 
& e^{S_1}e^{S_2}\\
& \left(\sqrt{\frac{a}{b}}-\sqrt{\frac{b}{a}}\right)e^{-S_1}+e^{-S_1}e^{-S_2} +e^{-S_1}e^{S_2} +\left(\sqrt{\frac{1}{a b}}+\sqrt{a b}\right) e^{S_2}\\
\end{array}\right),\nn
\eea
\bea
\label{rep4}
Y\to\left(\begin{array}{cc} 
\left(1+\frac{c d}{q}\right)\sqrt{a b}\, e^{S_1}- \sqrt{a b } \,c\, e^{-S_2}-  \sqrt{a b } \,d\, e^{2S_1}e^{S_2}& \\
-\frac{(1+a b)(c d+q)}{q}+ (1+a b)c\, e^{-S_1}e^{-S_2} + (1+a b)d\, e^{S_1}e^{S_2} +\sqrt{a b}\frac{q+c d}{q}e^{-S_1}
-\sqrt{a b}\,d\,e^{S_2}& \\
\end{array}\right.\nn \\
\left.\qquad\qquad \qquad\qquad\qquad \qquad\qquad\qquad \qquad
\begin{array}{cc} 
 &- \sqrt{a b } \,d\, e^{2S_1}e^{S_2}\\
& (1+a b)   d\, e^{S_1}e^{S_2} -\sqrt{a b}\, d\, e^{S_2}      \\
\end{array}\right), \nn
\eea
and
\bea
\label{rep5}
Z\to\left(\begin{array}{cc} 
\left(1+\frac{1}{a b}\right)\frac{q}{c} e^{S_1}e^{S_2}-\frac{1}{\sqrt{a b}\,c} e^{S_2} &\\
\frac{(1+a b)(c d+q)}{a b c d}-\left(1+\frac{1}{a b}\right)\left(\frac{q}{d} e^{-S_1}e^{-S_2}+
\frac{q}{c} e^{S_1}e^{S_2}\right)
-\frac{(c d+q)}{\sqrt{a b} \,c d} e^{-S_1} +
\frac{1}{\sqrt{a b}\,c} e^{S_2}&\\
\end{array}\right.\nn\\
\left.\qquad\qquad \qquad\qquad\qquad \qquad\qquad\qquad \qquad
\begin{array}{cc} 
&\frac{q^2}{\sqrt{a b} \, c} e^{2 S_1}e^{S_2}\\
 &\frac{(c d+q)}{\sqrt{a b} \,c d} e^{S_1}-
\frac{1}{\sqrt{a b}\,d} e^{-S_2}-\frac{q^2}{\sqrt{a b}\,c} e^{2 S_1}e^{S_2}
\\
\end{array}\right),\nn
\eea
where $S_1,S_2$ satisfy the following commutation relations:
\be\label{u0}
[S_{1},S_{2}]=i\pi \hbar,
\ee
and $q=e^{-i\pi \hbar}$,
gives and embedding of $\mathcal H$ into $Mat(2,\mathbb T_q)$. In particular, the images of $X,W,T,Z,T$ in  $GL(2,\mathbb T_q)$ satisfy the relations (\ref{LD0},\dots,\ref{LD4}) in the quantum ordering is dictated by the matrix product ordering\footnote{By this we mean that the product $A B$ of two matrices $A, B$ whose entries are in $\mathbb T_q$ is computed by keeping the entries of $A$ on the left matrix of the entries of $B$.}.
\end{theorem}

\proof First of all, note that since where $S_1,S_2$ satisfy the commutation relations (\ref{u0}), their exponentials satisfy the quantum torus commutation relations:
$$
e^{S_{2}} e^{S_{1}}= q e^{S_{1}} e^{S_{2}}, 
$$
for $q=e^{-i\pi \hbar}$ and therefore the images of $X,W,T,Z,T$ belong to  $Mat(2,\mathbb T_q)$. To prove that these images satisfy the relations(\ref{LD0},\dots,\ref{LD4}) 
 in which the quantum ordering is dictated by the matrix product ordering is a straightforward computation.

To prove that the map $H \to Mat(2,\mathbb T_q)$ defined by  (\ref{rep1}), (\ref{rep2}), (\ref{rep3}), (\ref{rep4}),(\ref{rep5}) is injective we need to prove that the images of 
$$
\left\{X^m Y^n\right\}_{n,m\in\mathbb Z}\cup\left\{T X^mY^n\right\}_{n,m\in\mathbb Z},
$$
are all linearly independent. To this aim, it is straightforward to prove by induction that for any $m\geq 0,n\geq 0$, $X^mY^n$ always contains terms with  
$$
e^{-m S_1},\dots, e^{- S_1},e^{ S_1}\dots, e^{(m+2 n-2) S_1},e^{-(m+n)S_2},\dots e^{- S_2},  e^{S_2},\dots e^{ (m+n-1) S_2},
$$
and  that for any $m\geq 0,n\geq 0$, $T X^mY^n$ always contains terms with  
$$
e^{-(m+1) S_1},\dots, e^{- S_1},e^{ S_1}, \dots, e^{(m+2 n-1) S_1},e^{-(m+n)S_2},\dots e^{- S_2},  e^{S_2},\dots,e^{ (m+n-1) S_2}.
$$
Since  $\{e^{k S_1},e^{m S_2}\}_{n,m\in\mathbb Z}$ are linearly independent, it automatically follows that the images of $\left\{X^m Y^n\right\}_{n,m\in\mathbb Z_{\geq 0}}$ and $\left\{T X^m Y^n\right\}_{n,m\in\mathbb Z_{\geq 0}}$ are all linearly independent. In a similar way it can be proved that the images of $\left\{T X^m Y^n\right\}_{n,m\in\mathbb Z_{\leq 0}}$ and 
$\left\{X^m Y^n\right\}_{n,m\in\mathbb Z_{\leq 0}}$ are all linearly independent. 
  
To show that all the images  of
$\left\{X^m Y^n\right\}_{n,m\in\mathbb Z}\cup \left\{TX^m Y^n\right\}_{n,m\in\mathbb Z}$ are linearly independent we proceed by contradiction. Assume that there exists a finite linear combination which gives zero: 
$$
\sum_{n,m} a_{n,m} X^m Y^n +\sum_{n,m} b_{n,m}  T X^m Y^n =0,
$$ 
take
$$
n_0=\min\{n| a_{n,m} \hbox{ or } b_{n,m}\neq 0\}, \quad m_0=\max\{l | a{n,m} \hbox{ or } b_{n,m} \neq 0\}.
$$
Iand multiply the above relation by $Y^{-n_0}X^{-m_0}$. Then we obtain a zero linear combination the set $\left\{T X^m  Y^n\right\}_{n,m\in\mathbb Z_{\geq 0}}\cup \left\{X^m  Y^n\right\}_{n,m\in\mathbb Z_{\geq 0}}$, which is absurd.
 \endproof

\begin{theorem}\label{th-rep-qT-PV}
The map given by the same formulae as Theorem \ref{th-rep-qT} for $T$, $X$ and $W$ and by
\be
Y\to \label{rep-lim1}
 \left(\begin{array}{cc}
 \sqrt{a b} \,e^{S_1}-\sqrt{a b}\,c\, e^{-S_2}&0\\
-(1+a b)+(1+a b)c\, e^{-S_1}e^{-S_2}+\sqrt{a b}\, e^{-S_1}&0\\
\end{array}\right)
\ee
\be
\label{rep-lim2}
Z\to \left(\begin{array}{cc}
0&0\\
\frac{a b+1}{a b} - c \frac{a b+1}{a b}e^{-S_1}e^{-S_2}-\frac{1}{\sqrt{a b}} e^{-S_1}&\frac{1}{\sqrt{a b}} e^{S_1}-\frac{c}{q\sqrt{a b}} e^{-S_2}
\end{array}\right)
\ee
where $S_1,S_2$ satisfy the commutation relations (\ref{u0}), gives and embedding of $\mathcal H_{V}$ into $Mat(2,\mathbb T_q)$. 
\end{theorem}
 
\proof We can prove that the images of $\left\{T X^m Y^n\right\}_{n\in\mathbb Z_{\geq 0},m\in\mathbb Z}\cup\left\{X^m Y^n\right\}_{n\in\mathbb Z_{\geq 0},m\in\mathbb Z}$ are all linearly independent in the same way as in the  proof of Theorem \ref{th-rep-qT}. To prove that the images of 
$$
\left\{T X^m Y^n\right\}_{n\in\mathbb Z_{\geq 0},m\in\mathbb Z}\cup\left\{X^m Y^n\right\}_{n\in\mathbb Z_{\geq 0},m\in\mathbb Z}\cup\left\{T X^m Z^n\right\}_{n\in\mathbb Z_{\geq 0},m\in\mathbb Z}\cup\left\{X^m Z^n\right\}_{n\in\mathbb Z_{\geq 0},m\in\mathbb Z}
$$
we proceed by contradiction: we assume that 
there exists a finite linear combination which gives zero: 
$$
\sum_{n,m} a_{n,m} X^m Y^n +\sum_{n,m} b_{n,m}  T X^m Y^n+\sum_{n,m} c_{n,m} X^m Z^n +\sum_{k,l} d_{n,m}  T X^m Z^n=0,
$$ 
in which at least one $c_{n,m}$ or $d_{n,m}$ is non zero for $n>0$.
Since the second column of $Y$ is identically zero, all elements in $\left\{T X^m Y^n\right\}_{n\in\mathbb Z_{> 0},m\in\mathbb Z}\cup\left\{X^m Y^n\right\}_{n\in\mathbb Z_{> 0},m\in\mathbb Z}$ have the second column identically equal to zero.  This means that the $12$ and $22$ elements of 
$$
\sum_{m} a_{0,m} X^m  +\sum_{m} b_{0,m}  T X^m +\sum_{n,m} c_{n,m} X^m Z^n +\sum_{k,l} d_{n,m}  T X^m Z^n
$$ 
must be $0$. We can prove by induction that the $12$ elements of  $X^m$ and $T X^m$ contain terms
$$
e^{- S_1},e^{ S_1}\dots, e^{m S_1}, e^{- S_2},  e^{S_2},\dots e^{m S_2},
$$
while the $22$ elements don't contain $e^{- S_2}$ either.
The $12$ and $22$ elements of  $X^m Z^n$ and $T X^m Z^n$ contain terms
$$
e^{- S_1},e^{ S_1}\dots, e^{(m+n) S_1}, e^{- S_2},  e^{S_2},\dots e^{m S_2}.
$$
By using the fact that  $\{e^{k S_1},e^{m S_2}\}_{n,m\in\mathbb Z}$ are linearly independent, we conclude. \endproof

\begin{theorem}\label{th-rep-qT-PIV}
The map:
\be
\label{rep-piv1}
T\to\left(\begin{array}{cc}
0&i b\, e^{S_1}\\
-i b\, e^{-S_1}&-1+b^2  \\
\end{array}\right),
\ee
\be
\label{rep-piv2}
X\to \left(\begin{array}{cc} - i\, e^{-S_1}+\frac{1}{q} e^{-S_1}e^{S_2} -i\left(b-\frac{1}{b}\right) e^{S_2}&-q\,e^{S_1}e^{S_2}\\
 i\, e^{-S_1}-\frac{1}{q} e^{-S_1}e^{S_2} +i\left(b-\frac{1}{b}\right) e^{S_2}&q\,e^{S_1}e^{S_2}\\
\end{array}\right)
\ee
\be
\label{rep-piv3}
W\to \left(\begin{array}{cc}e^{S_1}e^{S_2} &e^{S_1}e^{S_2} \\
 -i e^{-S_1}+e^{-S_1}e^{S_2} -i\left(b-\frac{1}{b}\right) e^{S_2}& -i e^{-S_1}+e^{-S_1}e^{S_2} -i\left(b-\frac{1}{b}\right) e^{S_2}\\
\end{array}\right)
\ee
\be
\label{rep-piv4}
Y\to \left(\begin{array}{cc} i b\, e^{S_1}-i b c\, e^{-S_2}&0\\
b^2-1+ i b\, e^{-S_1}+c(1-b^2) e^{-S_1}e^{-S_2} &0\\
\end{array}\right),
\ee
\be
\label{rep-piv5}
Z\to \left(\begin{array}{cc} 
0&0\\
\left(1-\frac{1}{b^2}\right)\left(1-c\, e^{-S_1}e^{-S_2} \right)+\frac{i}{b} e^{-S_1} &- \frac{i}{b} e^{S_1} +\frac{i c}{q b} e^{-S_2}\\
\end{array}\right),
\ee
where $S_1,S_2$ satisfy the commutation relations (\ref{u0}), gives and embedding  of $\mathcal H_{IV}$ into $Mat(2,\mathbb T_q)$. 
\end{theorem}
 
\proof The proof of this Theorem is very similar to the proof of Theorem \ref{th-rep-qT-PV}, except that in this case we need to prove that the images of 
\bea
&&
\left\{X^m Y^n\right\}_{m,n\in\mathbb Z_{\geq0}}\cup\left\{X^mT Y^n\right\}_{m,n\in\mathbb Z_{\geq0}}\cup 
\left\{X^m Z^n\right\}_{m,n\in\mathbb Z_{\geq0}}\cup\left\{X^mT Z^n\right\}_{m,n\in\mathbb Z_{\geq0}}\cup\nn\\
&&\cup\left\{W^m Y^n\right\}_{m,n\in\mathbb Z_{\geq0}}\cup\left\{W^mT Y^n\right\}_{m,n\in\mathbb Z_{\geq0}}\cup 
\left\{W^m Z^n\right\}_{m,n\in\mathbb Z_{\geq0}}\cup\left\{W^mT Z^n\right\}_{m,n\in\mathbb Z_{\geq0}}\nn
\eea
are all linearly independent. The only novelty is that now instead of negative powers of  $X$, we have positive powers of $W$. \endproof

\begin{theorem}\label{th-rep-qT-PIII}
The map given by the same formulae as Theorem \ref{th-rep-qT} for $T$, $X$ and $W$ and by
\be
Y\to \label{rep-lim1-III}
 \left(\begin{array}{cc}
\sqrt{a b}\, e^{-S_2}&0\\
-(1+a b) e^{-S_1}e^{-S_2}&0\\
\end{array}\right)
\ee
\be
\label{rep-lim2}
Z\to \left(\begin{array}{cc}
0&0\\
 \frac{a b+1}{a b}e^{-S_1}e^{-S_2}&\frac{1}{q\sqrt{a b}} e^{-S_2}
\end{array}\right)
\ee
where $S_1,S_2$ satisfy the commutation relations (\ref{u0}), gives and embedding of $\mathcal H_{III}$ into $Mat(2,\mathbb T_q)$. 
\end{theorem}
 
\proof The proof of this Theorem follows the same lines as the one of Theorem \ref{th-rep-qT-PV}. \endproof

\begin{theorem}\label{th-rep-qT-PII}
The map  given by the same formulae as Theorem \ref{th-rep-qT-PIV} for $T$, $Y$ and $Z$ and by
\be
\label{rep-pii2}
X\to \left(\begin{array}{cc}  \frac{1}{q} e^{-S_1}e^{S_2} -i\left(b-\frac{1}{b}\right) e^{S_2}&-q\,e^{S_1}e^{S_2}\\
 -\frac{1}{q} e^{-S_1}e^{S_2} +i\left(b-\frac{1}{b}\right) e^{S_2}&q\,e^{S_1}e^{S_2}\\
\end{array}\right)
\ee
\be
\label{rep-pii3}
W\to \left(\begin{array}{cc}e^{S_1}e^{S_2} &e^{S_1}e^{S_2} \\
e^{-S_1}e^{S_2} -i\left(b-\frac{1}{b}\right) e^{S_2}& e^{-S_1}e^{S_2} -i\left(b-\frac{1}{b}\right) e^{S_2}\\
\end{array}\right)
\ee
gives and embedding of $\mathcal H_{II}$ into $Mat(2,\mathbb T_q)$. 
\end{theorem}
 
\proof The proof of this Theorem follows closely the proof of Theorem \ref{th-rep-qT-PIV}, so we omit it. \endproof

\section{Confluent Zhedanov algebras and q-Askey scheme}\label{se:zhedanov-others}

In this Section we prove that the spherical sub-algebra of each confluent Cherednik algebra is isomorphic to the corresponding confluent Zhedanov algebra quotiented by its Casimir. Moreover we give a faithful representation of the confluent Zhedanov algebras and show that  they act as symmetries of some elements of the q-Askey scheme. Throughout this section many results on basic hypergeometric polynomials are used, they can be found in \cite{KLS} (see also \cite{AW} and \cite{GR} and references therein).
Throughout this Section we assume $q^m\neq 1,\,\forall m\in\mathbb Z$.

Before stating the isomorphism theorem we need to give the Casimir element for each confluent Zhedanov algebra:

\begin{lm}
For each index  $d= V,IV,III, III^{D_7}, III^{D_8}, II,I$, the Zhedanov algebra $\mathcal Z_d$  can be equivalently described as the algebra with two generators $K_0,K_1$ and two relations:
\bea\label{zhe-equiv}
&&
(q+q^{-1})K_1K_0K_1-K_1^2 K_0-K_0 K_1^2= B K_1 + \left(q-q^{-1}\right)^2 K_0 +D_0,\nn \\
&&
(q+q^{-1})K_0K_1K_0-K_0^2 K_1-K_1 K_0^2= B K_0 +D_1,
\eea
where the parameters $B$, $C_0$, $D_0$ and $D_1$ are chosen like in (\ref{eq:zhe-param}), 
and admits the following Casimir:
\bea
&&
\mathcal C=(K_1 K_0)^2-(q^2+1+q^{-2}) K_0 K_1 K_0 K_1+(q+q^{-1}) \left(q-{q^{-1}}\right)^2 K_0^2+\nn\\
&&
\qquad +(q+q^{-1} )K_0^2 K_1^2+B\left((q+1+q^{-1})K_0 K_1+K_1 K_0 \right)+\nn\\
&&
\qquad + (q+1+q^{-1})(D_0 K_0+D_1K_1).\nn
\eea
\end{lm}

Let us consider the quotiented Zhedanov algebra $\mathcal Z_d\slash\langle \mathcal C=\mathcal C_0\rangle$. Then we have the following:

\begin{theorem}\label{th:zhe-sph}
For each index  $d= V,IV,III, III^{D_7}, III^{D_8}, II,I$ the map: 
$$
i :\mathcal Z_d\slash\langle \mathcal C=\mathcal C_0\rangle \to  e\mathcal H_de,
$$
defined by
\bea\label{eq:isomor}
&&
i(K_0):= i\sqrt{a b} \hat X_2,\qquad i(K_1):= \hat X_1, \qquad i(1):= e\\
&&
i(K_2)=  i\sqrt{a b} \left(q-\frac{1}{q}\right) \hat X_3+
\frac{\sqrt{q}}{1+q} B \,e,\nn
\eea
where the parameters $B,C_0,D_0,D_1$ are given in terms of  the confluent Cherednik algebra parameters $a,b,c$ according to the formulae in table \ref{tab:sing} below,
is an algebra isomorphism.
\begin{table}[h]
\begin{center} 
\begin{tabular}{|c||c|c|c|c|} \hline 
 &B & $C_0$ & $D_0$& $D_1$\\ \hline 
$\mathcal Z_{V}$ & $\frac{(q-1)^2(c(1+a b/q)+a+b)}{q}$ & $\left(q-\frac{1}{q}\right)^2$&$- \frac{(q-1)^2(q+1)((a+b) c+q+a b)}{q^2}$&$ -\frac{(q-1)^2(q+1) a b c}{q^2}$ \\ \hline 
$\mathcal Z_{IV}$ &$\frac{(q-1)^2(c(q-b^2)+b q)}{q^2}$ &$0$ & $- \frac{(q-1)^2(q+1)}{q^2}b c$&$ \frac{(q-1)^2(q+1)}{q^2}b^2 c$ \\ \hline 
$\mathcal Z_{III}$  & $- \frac{(q-1)^2}{q}(1+a b/q)$ &  $\left(q-\frac{1}{q}\right)^2$&$\frac{(q-1)^2(q+1)}{q^2}(a+b)$&$0$ \\ \hline 
$\mathcal Z_{III}^{D_7}$ & $- \frac{(q-1)^2}{q}$ & $\left(q-\frac{1}{q}\right)^2$& $\frac{(q-1)^2(q+1)}{q^2}a$&$0$\\ \hline  
$\mathcal Z_{III}^{D_8}$  & $- \frac{(q-1)^2}{q}$ & $\left(q-\frac{1}{q}\right)^2$ & $0$ &$0$\\ \hline 
$\mathcal Z_{II}$  & $- \frac{(q-1)^2}{q^2}(q-b^2)\sqrt{q}$ &$0$ &$0$&$-\frac{(q-1)^2(q+1)}{q^2}b^2 \sqrt{q}$\\ \hline  
$\mathcal Z_{I}$ &$- \frac{(q-1)^2}{\sqrt{q}}$     &  $0$ &$0$&$0$\\ \hline
\end{tabular}
\vspace{0.2cm}
\end{center}
\caption{Values of the parameters $B,C_0,D_0,D_1$ in each confluent Zhedanov algebra in terms of  the confluent Cherednik algebra parameters $a,b,c$.}
\label{tab:sing}
\end{table}
\end{theorem}

\proof
To prove that $i$ is surjective it is enough to show that the relations (\ref{zhe-equiv}) for each $\mathcal Z_d$ are mapped by $i$ to the corresponding quantum commutation relations for $\hat X_1,\hat X_2,\hat X_3$ in $e\mathcal H_de$, where $d=II, III, III^{D_7}, III^{D_8}, IV,V$. 
This is a straightforward but lengthy computation that we omit for brevity. Note that in the case of $ Z_I(\mathcal C_0)$ the proof remains conjectural because both $a$ and 
 $b$ are zero in that case, so that we have no way to deduce $\hat X_3$ from $K_0,K_1,K_2$.
 
To prove that the map $i$ is injective, we only need to prove that in each spherical sub-algebra $e\mathcal H_d e$, where $d=II, III, III^{D_7}, III^{D_8}, IV,V$, the are no other relations a part from the quantum commutation relations for $\hat X_1,\hat X_2,\hat X_3$ and the quantum cubic relation. Let us illustrate how to proceed in the case of $e\mathcal H_V e$. We can use the cubic relation (\ref{eq:cubic-hatV}) to define $e$ in terms of $\hat X_1, \hat X_2, \hat X_3$. Then we use the quantum commutation relations (\ref{eq:skein-hat-5}) to order all words in the form $ \hat X_1^n \hat X_2^m \hat X_3^k$. If there is a further relation, this will necessarily have the following form:
$$
\sum_{m,n,k\in\mathbb Z_{\geq 0}} a_{m,n,k} \hat X_1^n \hat X_2^m \hat X_3^k =0.
$$
To prove that this can only be satisfied by choosing $a_{m,n,k} =0$ $\forall {m,n,k\in\mathbb Z}$, we use the embedding
Theorem \ref{th-rep-qT-PV} and show by induction the image of each element $ \hat X_1^n \hat X_2^m \hat X_3^k $ contains the following terms:
$$
e^{-(k+n+1) S_1},\dots, e^{-S_1},e^{S_1},\dots, e^{(n+m+1) S_1},
e^{-(2 k+n+m) S_2},\dots, e^{-S_2},e^{S_2},\dots, e^{n S_2}.
$$
Since $\{e^{k S_1},e^{m S_2}\}_{n,m\in\mathbb Z}$ are linearly independent, it  follows that the images of $ \hat X_1^n \hat X_2^m \hat X_3^k $ are all linearly independent. In fact let $n_0=\max\{n: a_{m,n,k}\neq 0\}$. If $n_0>0$, then we have terms $e^{n_0 S_1}$ that can't be eliminated by any other term. This proves that $n_0=0$. Let $m_0=\max\{m: a_{m,n,k}\neq 0\}$, then $m_0=0$ as otherwise we have terms with $e^{(m_0+1) S_2}$ that can't be eliminated by any other term. Similarly,  $k_0=\max\{k: a_{m,n,k}\neq 0\}$, then $k_0=0$ as otherwise we have terms with $e^{-(k_0+1) S_1}$ that can't be eliminated by any other term.

The proof for all other cases   $e\mathcal H_{IV} e$,  $e\mathcal H_{III}e$, $e\mathcal H_{II} e$ follows the same reasoning and we omit it.

To prove the statement in the case of  the two algebras $e\mathcal H_{III^{D_7}}e$ and $e\mathcal H_{III^{D_8}}e$, we observe that they are special cases of $e\mathcal H_{III}e$ for $b=0$ and $a=b=0$ respectively.  Bearing this in mind, we can prove that  the two algebras $e\mathcal H_{III^{D_7}}e$ and $e\mathcal H_{III^{D_8}}e$ can be embedded in to 
$Mat_2(\mathbb T_q)$ and follow the same reasoning as above.
 \endproof

Now in each case we give a faithful representation the confluent Zhedanov algebras either on the space of symmetric Laurent polynomials $\mathcal L_{sym}$ or on the space of polynomials $\mathcal P$. In order to prove that our representation is faithful,  we need first  the following lemma:

\begin{lm}\label{lm:basisV}
The quotiented Zhedanov algebra $\mathcal Z_d\slash\langle \mathcal C=\mathcal C_0\rangle$, has elements
$$
K_0^n(K_1K_0)^lK_1^m, \qquad m,n= 0,1, 2, 3,\dots, \quad l=0,1,
$$
as a basis.
\end{lm}

\proof The proof is the same as in the non--confluent case see  the results contained in Section 2 of \cite{K1}.
\endproof

We will give the proofs only in the first case, i.e. for  $\mathcal Z_V$. All other proofs are similar and therefore will be omitted for brevity.

\subsection{Representation of $\mathcal Z_V$ and continuous dual q-Hahn polynomials}

\begin{lm}
The confluent Zhedanov algebra  $\mathcal Z_V\slash\langle \mathcal C=\mathcal C_0\rangle$ admits the following representation on the space $\mathcal L_{sym}$ of symmetric Laurent polynomials:
\bea
\label{eq:basicV1}
&\qquad (K_1 f)[x]:=&( x+\frac{1}{x}) f[x],\\
\label{eq:basicV2}
&\qquad(K_0 f)[x]:=& \frac{(1-a x)(1-b x)(1-c x)}{(1-x^2)(1-q x^2)}(f[qx]-f[x])+ f[x]-\nn\\
&&-x\frac{(a-x)(b-x)(c-x)}{(1-x^2)(q-x^2)}(f[q^{-1}x]-f[x]).
\eea
\end{lm}

\proof By expressing the confluent Zhedanov algebra structure constants by the parameters $a,b,c$:
it is a straightforward computation (see notebook 8 in \cite{MN}) to prove that the operators satisfy the relations (\ref{zhe1-pv},\ref{zhe2-pv},\ref{zhe3-pv}),with constants $B,C_1,D_0,D_1$ given in table 1. 
\endproof

\begin{lm}
The continuous dual $q$-Hahn  polynomials:
$$
p_n(x;a,b,c,d) := \frac{(ab,ac;q)_n}{a^n}    {}_3\phi_2\left(\begin{array}{cc}
q^{-n},a x,a x^{-1}\\
ab,ac,
\end{array};q,q
\right),
$$
where $a=-\frac{u_1}{k_1}$, $b= u_1 k_1$, $c=-\frac{\sqrt{q}}{u_0}$, are eigenfunctions of the $K_0$ operator:
$$
K_0 p_n =q^{-n}  p_n.
$$
\end{lm}

\proof Note that the confluent Zhedanov algebra  $\mathcal Z_V\slash\langle \mathcal C= \mathcal C_0\rangle$ is obtined as the limit for $d\to 0$ of the general Zhedanov algebra   $\mathcal Z\slash\langle  \mathcal C= \mathcal C_0\rangle$. Analogously, the representation (\ref{eq:basicV1}, \ref{eq:basicV2}) is the limit as $d\to 0$ of the representation (\ref{eq:basic4},\ref{eq:basic5}) of the general Zhedanov algebra and the  continuous dual $q$-Hahn  polynomials are obtained as limit for $d\to 0$ of the Askey--Wilson polynomilas. As a consequence this proof follows from the fact that the Askey--Wilson polynomials are eigenfunctions of the operator $K_0$ defined by (\ref{eq:basic5}).
\endproof

\begin{lm}
The representation (\ref{eq:basicV1},\ref{eq:basicV2}) is faithful.
\end{lm}

\proof
This proof follows the same lines as the proof of Theorem 2.2 in \cite{K1}, where we replace the Askey-Wilson polynomials by the dual $q$-Hahn polynomials.
\endproof

\subsection{Big $q$--Jacobi polynomials}

In order to understand the relation between the confluent Zhedanov algebra $\mathcal Z_V$ and the Big q-Jacobi polynomials, we use the confluent version of the Cherednik algebra automorphism $\gamma$ defined in sub--section \ref{suse:aut}. As shown in Section \ref{der}, in the confluent limit, $\gamma$ defines an algebra isomorphism (\ref{eq:gamma-V}) mapping $\mathcal H_V$  to $\mathcal H_{V}^\gamma$ given by (\ref{sahi1-Vgamma},\ref{sahi2-Vgamma}\ref{sahi3-Vgamma}\ref{sahi4-Vgamma}\ref{sahi6-Vgamma}). On the confluent Zhedanov algebra we obtain the following result:

\begin{lm}\label{auto-Zhe}
The transformation
$$
\gamma(K_0,K_1)=\left(k_1 u_1 K_1,\frac{1}{u_1}\left(K_0+\frac{q^{\frac{3}{2}}}{(q+1)(q-1)^2}[K_2,K_1]\right)
\right),
$$
is an isomorphism mapping $\mathcal Z_V$ to  $\mathcal Z_V^\gamma$, which is the algebra generated by $K_0^\gamma,K_1^\gamma$ with relations:
\bea
&&\label{eq:zhe-parV-gamma1}
(q+q^{-1})K_1^\gamma K_0^\gamma K_1^\gamma -(K_1^\gamma)^2 K_0^\gamma-K_0^\gamma (K_1^\gamma)^2=
 B^\gamma K_1^\gamma +D_0^\gamma, \\
&&
(q+q^{-1})K_0^\gamma K_1^\gamma K_0^\gamma -(K_0^\gamma)^2 K_1^\gamma-K_1^\gamma (K_0^\gamma)^2= 
B^\gamma K_0^\gamma +C_1^\gamma K_1^\gamma+D_1^\gamma,\nn
\eea
for some arbitrary constants $B^\gamma, C_1^\gamma, D_0^\gamma $ and $D_1^\gamma$\footnote{Only three of these are independent, as one constant can be rescaled to any arbitrary value by rescaling the generators.}.\end{lm}

\proof
It is a straightforward computation to prove that on the spherical sub--algebra $e\mathcal H_{V}e$ the isomorphism $\gamma$ acts as follows:
$$
\gamma( \hat X_1,\hat X_2,\hat X_3)=\left( \frac{\sqrt{q}}{q-1}[\hat X_3,\hat X_1]+\hat X_2,\hat X_1,\hat X_3\right).
$$
By using the isomorphism
$i :
\mathcal Z_V\slash\langle \mathcal C=\mathcal C_0\rangle
 \to  e\mathcal H_V e$ obtained in theorem \ref{th:zhe-sph} we can deduce the result.
\endproof

\begin{lm}
The confluent Zhedanov algebra  $\mathcal Z_V^\gamma\slash\langle \mathcal C=\mathcal C_0\rangle$ admits the following representation on the space $\mathcal P$ of  polynomials:
\bea
\label{eq:basicV1-J}
&\qquad (K_1^\gamma f)[x]:=& x \, f[x],\\
\label{eq:basicV2-J}
&\qquad(K_0^\gamma f)[x]:=&\frac{q\left (\lambda \tilde c x+ \tilde a(x(1+ \tilde b)- \tilde c(1+q-\lambda x))\right)}{\lambda^2 x^2} f[x]+\\
&&+\frac{(\lambda x-q \tilde  a)(\lambda x-q  \tilde c)}{\lambda^2 x^2} f\left[\frac{x}{q}\right]+\frac{q(\lambda x-1) \tilde a( \tilde b\lambda x - \tilde c)}{x^2}f[q x].\nn
\eea
\end{lm}

\proof Indeed the generators defined by (\ref{eq:basicV1-J}, \ref{eq:basicV2-J}) satisfy the relations (\ref{eq:zhe-parV-gamma1}) for 
$$
B^\gamma = (q-1)^2\frac{ \tilde c+ \tilde a(1+ \tilde b+ \tilde c)}{\lambda},\qquad D_0^\gamma=-(q+1)(q-1)^2\frac{  \tilde a  \tilde c}{\lambda^2},
$$
$$
C_1^\gamma= q\left(q-\frac{1}{q}\right)^2 \tilde a  \tilde b,\qquad 
D_1^\gamma= -(q-1)^2(q+1) \frac{ \tilde a ( \tilde c+ \tilde b(1+ \tilde a+ \tilde  c))}{\lambda}.
$$
see notebook 9 in \cite{MN}).
\endproof
 
The proof of the following two results is obtained by taking substituting $x\to \frac{x}{\eps}$, $a\to\eps\lambda$, $b\to\frac{\tilde a q}{\eps \lambda}$, $c\to\frac{\tilde c q}{\eps \lambda}$, $d\to\eps \lambda \frac{\tilde b}{\tilde c}$ and taking the limit as $\eps\to 0$ in the analogous results for the Askey Wilson polynomials in \cite{K1}.

\begin{lm}
The big $q$--Jacobi polynomials:
$$
P_n(x;a,b,c,d) :=   {}_3\phi_2\left(\begin{array}{cc}
q^{-n},a b q^{n+1}, x\\
a q, c q,
\end{array};q,q
\right),
$$
are eigenfunctions of the $K_0$ operator:
$$
K_0 P_n[\lambda x] =\frac{1+q^{2n+1}a b}{q^n}  P_n[\lambda x].
$$
\end{lm}

\begin{lm}
The representation (\ref{eq:basicV1-J},\ref{eq:basicV2-J}) is faithful.
\end{lm}

\begin{remark}
We have introduced an arbitrary parameter $\lambda$ in our representation to account for the freedom of direction of confluence. This freedom will be used later on to produce exact matches between parameters in the polynomials and parameters in the confluent Zhedanov algebras,
\end{remark}

\subsection{Representation of $\mathcal Z_{IV}$ and  Big $q$--Laguerre Polynomials}

Note that the algebra $\mathcal H_{IV}$ can be obtained as limit of the algebra  $\mathcal H_{V}^\gamma$  by rescaling 
$V_0^\gamma\to\frac{1}{\eps}V_0^\gamma$, $\check{V_0}^\gamma\to\frac{1}{\eps}\check{V_0}^\gamma$, $k_0\to\eps$, $u_0\to\eps u_0$. This shows that the confluent Zhedanov algebra  $\mathcal Z_{IV}$   can be obtained from  $\mathcal Z_{V}^\gamma$ in the limit $b\to 0$. This leads to the following results:

\begin{lm}
The confluent Zhedanov algebra  $\mathcal Z_{IV}\slash\langle \mathcal C=\mathcal C_0\rangle$ admits the following representation  on the space $\mathcal P$ of  polynomials:
\bea
\label{eq:basicIV1}
&\qquad (K_1 f)[x]:=& x \, f[x],\\
\label{eq:basicIV2}
&\qquad(K_0 f)[x]:=&\frac{q\left (\lambda\tilde c x+\tilde a(x -\tilde c(1+q-\lambda x))\right)}{\lambda^2 x^2} f[x]+\\
&&+\frac{(\lambda x-q \tilde a)(\lambda x-q \tilde c)}{\lambda^2 x^2} f\left[\frac{x}{q}\right]-\frac{q(\lambda x-1)\tilde a \tilde c}{x^2}f[q x].\nn
\eea
\end{lm}

\proof Indeed the generators defined by (\ref{eq:basicIV1}, \ref{eq:basicIV2}) satisfy the relations (\ref{zhe-equiv}) for 
$$
B = (q-1)^2\frac{\tilde c+\tilde a(1+\tilde c)}{\lambda},\qquad D_0^\gamma=-(q+1)(q-1)^2\frac{\tilde  a\tilde  c}{\lambda^2},
$$
$$
D_1 = -(q-1)^2(q+1) \frac{\tilde a \tilde c}{\lambda},
$$
(see notebook 10 in \cite{MN}). 
Choosing $\tilde a= -\frac{b^2}{q}$, $\tilde c=-\frac{b c}{q}$, $\lambda= -b$, we find that these formulae are the same as the ones in Table 1.
\endproof

\begin{lm}
The big $q$--Laguerre polynomials:
$$
P_n(x;\tilde a,\tilde c,q) :=   {}_3\phi_2\left(\begin{array}{cc}
q^{-n},0, x\\
\tilde a q, \tilde c q,
\end{array};q,q
\right),
$$
are eigenfunctions of the $K_0$ operator:
$$
K_0 P_n[\lambda x] =q^{-n}  P_n[\lambda x].
$$
\end{lm}

\begin{lm}
The representation (\ref{eq:basicIV1},\ref{eq:basicIV2}) is faithful.
\end{lm}

\subsection{Representation of $\mathcal Z_{III}$ and Al--Salam-Chihara Polynomials}

\begin{lm}
The confluent Zhedanov algebra  $\mathcal Z_{III}$ admits the following representation on the space $\mathcal L_{sym}$ of symmetric Laurent polynomials:
\bea
\label{eq:basicIII1}
&\qquad (K_1 f)[x]:=&( x+\frac{1}{x}) f[x],\\
\label{eq:basicIII2}
&\qquad(K_0 f)[x]:=&x  \frac{ (1-a x)(1-b x)}{(1-x^2)(1-q x^2)}(f[qx]-f[x])+\nn\\
& &+x \frac{(a-x)(b-x)}{(1-x^2)(q-x^2)}(f[q^{-1}x]-f[x]).
\eea
\end{lm}

\proof By expressing the confluent Zhedanov algebra structure constants by the parameters $a,b,c$ as in Table 1, it is a straightforward computation (see notebook 11 in \cite{MN}) to prove that the operators (\ref{eq:basicIII1},\ref{eq:basicIII2}) satisfy (\ref{zhe-equiv}).
\endproof

\begin{lm}
The Al-Salam-Chihara polynomials:
$$
Q_n(x;a,b,d) := \frac{(ab;q)_n}{a^n}    {}_3\phi_2\left(\begin{array}{cc}
q^{-n},a x,a x^{-1}\\
ab,0,
\end{array};q,q
\right),
$$
are eigenfunctions of the following operator
$$
K_0^\beta := \frac{q}{q^2-1}\left(K_0 K_1-q K_1K_0 -\frac{(a+b)(q-1)}{q}\right),
$$
with eigenvalues 
$$
\frac{1}{q^n}-1+ \frac{1+a+b-a b}{q+1}.
$$
\end{lm}

\begin{remark}
Note that now the operator $K_0$ does not act nicely on the Al-Salam-Chihara polynomials; we had to replace it by the new operator $K_0^\beta$. This is due to the fact that  in terms of generators $T_0,T_1,X,W$ and parameters $a,b$, the algebra $\mathcal H_{III}$ is obtained as a limiting case of $\mathcal H_{V}$ for $c\to\infty$, while the Al--Salam--Chihara polynomials are obtained by the continuous dual $q$--Hahn polynomials in the limit $c\to 0$. In order to match these two asymptotics, we need to perform the following transformation of the $\mathcal H_{III}$ algebra:
$$
\beta(T_0,T_1,X)= \left(\sqrt{q} W T_0-1,T_1,X\right).
$$
This transformation is  an isomorphism between the  algebra $\mathcal H_{III}$  in the representation (\ref{sahi1-III},\dots,\ref{sahi6-III}) and the algebra  $\mathcal H_{III}^\beta$ generated by $T_0,T_1,X,W$ and relations (\ref{sahi1-III}), (\ref{sahi2-III}), (\ref{sahi4-III}) and the algebra:
\bea
\label{sahi3-1-III}
T_0^2+T_0=0,\\
\label{sahi6-1-III}
q T_0 W= X (T_0+1),
\eea
which is obtained from the $\mathcal H_{V}$ algebra as $c\to 0$. This transformation acts on the confluent Zhedanov algebra $\mathcal Z_{III}$ as follows:
$$
\beta(K_0,K_1)= (K_0^\beta ,K_1).
$$
\end{remark}

\begin{lm}
The representation (\ref{eq:basicIII1},\ref{eq:basicIII2}) is faithful.
\end{lm}

\proof We use the same idea as in the remark above and prove all statements as limits for $c\to 0$ of the analogous statements for the algebra $\mathcal H_{V}$  and the dual q--Hahn polynomials.\endproof 

\subsection{Representation of $\mathcal Z_{III}^{D_7}$ and  continuous Big $q$-Hermite Polynomials}
Note that the algebra $\mathcal H_{III}^{D_7}$ can be obtained as limit of the algebra  $\mathcal H_{III}$  by taking the limit $b\to 0$ and $c\to 1$ (see notebook 12 in \cite{MN}). This leads to the following results:

\begin{lm}
The confluent Zhedanov algebra admits the following representation on the space $\mathcal L_{sym}$ of symmetric Laurent polynomials:
\bea
\label{eq:basic4D7}
&\qquad (K_1 f)[x]:=&( x+\frac{1}{x}) f[x],\\
\label{eq:basic5D7}
&\qquad(K_0 f)[x]:=& x\frac{ (1-a x)}{(1-x^2)(1-q x^2)}(f[qx]-f[x])-\nn\\
& &  x^2\frac{(a-x)}{(1-x^2)(q-x^2)}(f[q^{-1}x]-f[x]).
\eea
\end{lm}

\begin{lm} The continuous big $q$-Hermite polynomials:
$$
H_n(x;a,b,c,d) := \frac{1}{a^n}    {}_3\phi_2\left(\begin{array}{cc}
q^{-n},a x,a x^{-1}\\
0,0,
\end{array};q,q
\right),
$$
are eigenfunctions of the following operator
$$
K_0^\beta := \frac{q}{q^2-1}\left(K_0 K_1-q K_1K_0 -\frac{a (q-1)}{q}\right),
$$
with eigenvalues 
$$
\frac{1}{q^{n+1}}- \frac{1+a}{q+1}.
$$
\end{lm}

\begin{lm}
The representation (\ref{eq:basic4D7},\ref{eq:basic5D7}) is faithful.
\end{lm}

\subsection{Representation of $\mathcal Z_{III}^{D_8}$ and continuous $q$-Hermite Polynomials}

Note that the algebra $\mathcal H_{III}^{D_8}$ can be obtained as limit of the algebra  $\mathcal H_{III}^{D_7}$  by taking the limit $a\to 0$. This leads to the following results  (see notebook 13 in \cite{MN}):

\begin{lm}
The confluent Zhedanov algebra  $\mathcal Z_{III}^{D_8}$ admits the following representation on the space $\mathcal L_{sym}$ of symmetric Laurent polynomials:
\bea
\label{eq:basic4D8}
&\qquad (K_1 f)[x]:=&( x+\frac{1}{x}) f[x],\\
\label{eq:basic5D8}
&\qquad(K_0 f)[x]:=&  x\frac{ 1}{(1-x^2)(1-q x^2)}(f[qx]-f[x])+\nn\\
& &  x^3\frac{1}{(1-x^2)(q-x^2)}(f[q^{-1}x]-f[x]).
\eea
\end{lm}

\begin{lm} The continuous $q$-Hermite polynomials:
$$
H_n(x;a,b,c,d) := x^n   {}_2\phi_0\left(\begin{array}{cc}  
q^{-n},0\\
- \\
\end{array};q,\frac{q^n}{x^2}
\right),
$$
are eigenfunctions of the following operator
$$
K_0^\beta := \frac{q}{q^2-1}\left(K_0 K_1-q K_1K_0\right),
$$
with eigenvalues 
$$
\frac{1}{q^{n+1}}- \frac{1}{q+1}.
$$
\end{lm}

\begin{lm}
The representation (\ref{eq:basic4D8},\ref{eq:basic5D8}) is faithful.
\end{lm}

\subsection{Representation of $\mathcal Z_{II} $ and little $q$--Laguerre/Wall polynomials}

\begin{lm}
The confluent Zhedanov algebra  $\mathcal Z_{II}\slash\langle \mathcal C=\mathcal C_0\rangle$ admits the following representation on the space $\mathcal P$ of  polynomials:
\bea
\label{eq:basicII1}
&\qquad (K_1 f)[x]:=& x \, f[x],\\
\label{eq:basicII2}
&\qquad(K_0 f)[x]:=&\frac{1+\tilde a}{\lambda x} f[x]+\frac{\lambda x-1}{\lambda x} f\left[\frac{x}{q}\right]-\frac{\tilde a}{\lambda x} f[q x].
\eea
\end{lm}

\proof Indeed the generators defined by (\ref{eq:basicII1}, \ref{eq:basicII2}) satisfy the relations (\ref{zhe-equiv}) with
\be
B=\frac{(q-1)^2(1+\tilde a)}{\lambda q},\qquad D_1= -\frac{(q-1)^2(1+q)}{\lambda q} \tilde a.\nn
\ee
Choosing $\tilde a= -\frac{b^2}{q}$,  $\lambda= -b$, we find that these formulae are the same as the ones in Table 1.
\endproof

The following results can be proved by taking $c\to -\frac{1}{\eps}$ and $x\to \frac{q x}{\eps}$ and letting $\eps\to 0$ in the results proved for 
$\mathcal Z_{IV}$ (see also notebook 14 in \cite{MN}).

\begin{lm}
The little $q$--Laguerre polynomials:
$$
p_n(x;a,c,d) :=   {}_3\phi_2\left(\begin{array}{cc}
q^{-n},0\\
a q, 
\end{array};q,q x
\right),
$$
are eigenfunctions of the $K_0$ operator:
$$
K_0 p_n[ x] =q^{-n}  p_n[x].
$$
\end{lm}

\begin{lm}
The representation (\ref{eq:basicII1},\ref{eq:basicII2}) is faithful.
\end{lm}

\subsection{Representation of $\mathcal Z_{I} $ and a special case of the little $q$--Laguerre/Wall polynomials}

The confluent Zhedanov algebra  $\mathcal Z_{I}$  can be obtained from  $\mathcal Z_{II}$ in the limit $a\to 0$. This leads to the following results:

\begin{lm}
The confluent Zhedanov algebra  $\mathcal Z_{I}\slash\langle \mathcal C=\mathcal C_0\rangle$ admits the following representation on the space $\mathcal P$ of  polynomials:
\bea
\label{eq:basicI1}
&\qquad (K_1 f)[x]:=& x \, f[x],\\
\label{eq:basicI2}
&\qquad(K_0 f)[x]:=&\frac{1}{x} f[x]+\frac{x-1}{x} f\left[\frac{x}{q}\right],\nn
\eea
with
\be
B=\frac{(q-1)^2}{q}.
\ee
\end{lm}

\begin{lm}
The little $q$--Laguerre polynomials with $a=0$:
$$
p_n(x;0,c,d) :=   {}_3\phi_2\left(\begin{array}{cc}
q^{-n},0\\
0, 
\end{array};q,q x
\right),
$$
are eigenfunctions of the $K_0$ operator:
$$
K_0 p_n[ x] =q^{-n}  p_n[x].
$$
\end{lm}

\begin{lm}
The representation (\ref{eq:basicI1},\ref{eq:basicI2}) is faithful.
\end{lm}


\setcounter{section}{0}

\def\thetheorem{A.\arabic{theorem}}
\def\theprop{A.\arabic{prop}}
\def\thelemma{A.\arabic{lm}}
\def\thecor{A.\arabic{cor}}
\def\theexam{A.\arabic{exam}}
\def\theremark{A.\arabic{remark}}
\def\theequation{A.\arabic{equation}}

\appendix{The Painlev\'e monodromy manifolds}\label{se:sing}

All the Painlev\'e differential equations arise as monodromy preserving deformations of an auxiliary system of two first order ODEs. The monodromy data of such auxiliary system are encoded in their monodromy manifolds which can all be described by cubic surfaces in $\mathbb C^3$ defined by the zero locus of the following polynomials in
 $\mathbb C[x_1,x_2,x_3]$:
\begin{eqnarray}\nn
PVI & \begin{array}{l}
x_1x_2x_3 + x_1^2+x_2^2+x_3^2-(G_1 G_\infty+G_2 G_3) x_1 -(G_2 G_\infty+G_1 G_3) x_2  -\\
\qquad -(G_3 G_\infty+G_1 G_2) x_3+ G_1^2+G_2^2+G_3^2+G_\infty^2+G_1 G_2G_3 G_\infty-4,\\
\end{array}\nn\\
&&\nn\\
PV & \begin{array}{l}
x_1x_2x_3 + x_1^2+x_2^2-(G_1 G_\infty+G_2 G_3) x_1 -(G_2 G_\infty+G_1 G_3) x_2  -\\
\qquad -G_3 G_\infty x_3+ G_\infty^2+G_3^2+G_1 G_2G_3 G_\infty, \\
\end{array}\nn\\
&&\nn\\
PIV& \begin{array}{l}x_1x_2x_3 + x_1^2-(G_1 G_\infty+G_2 G_3) x_1 -G_2 G_\infty  x_2  -G_3 G_\infty x_3+\\
\qquad+ G_\infty^2+G_1 G_2G_3 G_\infty, \\
\end{array}
\nn\\
&&\nn\\
PIII &x_1 x_2 x_3 + x_1^2+x_2^2-G_1 G_\infty x_1 -G_2 G_\infty x_2 + G_\infty^2,\nn\\
&&\nn\\
PIII^{D_7} & x_1x_2x_3 + x_1^2+x_2^2-G_1 G_\infty x_1 -G_2 G_\infty x_2,\nn\\
&&\nn\\
PIII^{D_8} &x_1x_2x_3 + x_1^2+x_2^2 -G_2 G_\infty x_2,\nn\\
&&\nn\\
PII&x_1x_2x_3 + x_1^2-G_1 G_\infty x_1 -G_2 G_\infty  x_2  + G_\infty^2,\nn\\
&&\nn\\
PI& x_1x_2x_3 -G_1 G_\infty x_1 -G_2 G_\infty  x_2 + G_\infty^2,\nn
\end{eqnarray}
where $G_1,G_2,G_3,G_\infty$ are some constants related to the parameters appearing in the corresponding Painlev\'e equation as described in  \cite{CMR}. Note that not all parameters are independent is the above cubic polynomials: $G_\infty$ can be chosen arbitrarily in the $PV$, $PIII$, $PIII^{D_7}$, $PIII^{D_8}$, $PII$ and $PI$ cubics, $G_2$ can be chosen arbitrarily in the $PIII^{D_7}$, $PIV$, $PII$, $PI$ cubics, $G_3$ can be chosen arbitrarily in the $PII$,  and $PI$ cubics. To match these cubics to the ones found as semiclassical limits of the spherical sub-algebras of the confluent Cherednik algebras defined in this paper we need to perform the following substitutions:
$$
x_1=\left\{\begin{array}{lc}
-X_2&\hbox{for } PVI, PV, PIII, PIII^{D_7}, PIII^{D_8},\\
-X_3&\hbox{for } PIV,PII,PI,\\
\end{array}\right.
$$
$$
x_2=\left\{\begin{array}{lc}
-X_3&\hbox{for }  PVI, PV, PIII, PIII^{D_7}, PIII^{D_8},\\
-X_1&\hbox{for } PIV,PII,PI,\\
\end{array}\right.
$$
$$
x_3=\left\{\begin{array}{lc}
-X_1&\hbox{for }  PVI, PV, PIII, PIII^{D_7}, PIII^{D_8},\\
-X_2&\hbox{for } PIV,PII,PI,\\
\end{array}\right.
$$
$$
G_\infty=\left\{\begin{array}{lc}
-i \overline k_0&\hbox{for } PVI,\\
i&\hbox{for } PV, PIII, PIII^{D_7}, PIII^{D_8},PII\\
i c&\hbox{for } PIV,\\
1&\hbox{for } PI,\\
\end{array}\right.
$$
$$
G_3=\left\{\begin{array}{lc}
-i \overline u_0&\hbox{for } PVI,\\
-i c&\hbox{for } PV,\\
0&\hbox{for } PIII, PIII^{D_7}, PIII^{D_8},PII,PI\\
-i &\hbox{for } PIV,\\
\end{array}\right.
$$
$$
G_2=\left\{\begin{array}{lc}
-i\overline k_1&\hbox{for } PVI,  PV,\\
i\left(\sqrt{a b}+\frac{1}{\sqrt{a b}}\right) &\hbox{for } PIII,\\
i&\hbox{for } PIII^{D_7}, PIII^{D_8},PII,\\
- i &\hbox{for } PIV,\\
-1 &\hbox{for } PI,\\
\end{array}\right.
$$
$$
G_1=\left\{\begin{array}{lc}
-i\overline u_1&\hbox{for }PVI,  PV,\\
i\left(\sqrt{\frac{b}{a}}+\sqrt{\frac{a}{b}}\right) &\hbox{for } PIII,\\
-i\left(b-\frac{1}{b}\right) &\hbox{for } PIV,\\
i\left(b-\frac{1}{b}\right) &\hbox{for } PII,\\
-i a&\hbox{for } PIII^{D_7},\\
0 &\hbox{for }  PIII^{D_8},\\
1 &\hbox{for }  PI.\\
\end{array}\right.
$$

\vskip 2mm \noindent{\bf Acknowledgements.} The author is specially grateful to T. Koornwinder for reading trough this paper and coming up with many useful comments. The author wishes to thank also M. Balagovic, Yu. Berest, O. Chalyck, F. Eshmatov, P. Etingof, N. Joshi,  C. Korff , M. Noumi, V. Rubtsov, J. Stokman and P. Terwilliger  for interesting discussions on this subject.
This research was supported by the EPSRC and by the Hausdorff Institute. The author wishes to thank the referees for their useful comments.

\end{document}